\documentclass[aos]{imsart}

\RequirePackage{amsthm,amsmath,amsfonts,amssymb}
\RequirePackage[numbers]{natbib}
\RequirePackage[colorlinks,citecolor=blue,urlcolor=blue]{hyperref}
\RequirePackage{graphicx}
\RequirePackage{xcolor}
\RequirePackage{dsfont}
\RequirePackage{tikz,caption,subcaption,float,enumerate}
\RequirePackage{booktabs}
\RequirePackage{comment}

\startlocaldefs
\newcommand\norm[1]{\left\lVert#1\right\rVert}
\newcommand\set[1]{\{ \, #1 \,\}}

\theoremstyle{plain}
\newtheorem{theorem}{Theorem}[section]
\newtheorem{lemma}{Lemma}[section]
\newtheorem{remark}{Remark}[section]
\newtheorem{definition}{Definition}[section]

\renewcommand{\L}{\mathcal{L}}
\newcommand{\MS}{\mathcal{M}}%metric space
\newcommand{\D}{\mathcal{D}}
\newcommand{\R}{\mathbb{R}}

\newcommand{\Diag}{\Delta}
\endlocaldefs

\begin{document}

\begin{frontmatter}
\title{Kernel Regression using Persistence Diagrams}
%\title{A sample article title with some additional note\thanksref{T1}}
\runtitle{Kernel Regression using Persistence Diagrams}
%\thankstext{T1}{A sample of additional note to the title.}

\begin{aug}
%%%%%%%%%%%%%%%%%%%%%%%%%%%%%%%%%%%%%%%%%%%%%%%
%% Only one address is permitted per author. %%
%% Only division, organization and e-mail is %%
%% included in the address.                  %%
%% Additional information such as            %%
%% identifying the corresponding author must %%
%% be included in in the Acknowledgments     %%
%% section if necessary.                     %%

%I have just used the alphabetical ordering. Please feel free to change!
\author[A]{\fnms{Saksham}~\snm{Aggarwal}\ead[label=e1]{saksham94.sa@gmail.com}},
\author[A]{\fnms{James}~\snm{Nichols}\ead[label=e2]{james.nichols@anu.edu.au}},
\and
\author[A]{\fnms{Katharine}~\snm{Turner}\ead[label=e3]{katharine.turner@anu.edu.au}}

\address[A]{Mathematical Sciences Institute, Australian National University \printead[presep={,\ }]{e1,e2,e3}}
\end{aug}

\begin{abstract}
Persistence diagrams are widely used in topological data analysis to summarize topological and geometric structure in data. We study nonparametric regression where inputs are persistence diagrams and outputs are real-valued. Using kernel-based local averaging estimators defined directly over diagram space with Wasserstein distance, we establish consistency and derive bounds on the rate of convergence under suitable assumptions in the $L^2$-theory. We also present computational experiments that demonstrate strong performance of the proposed estimator, particularly in low data regimes.
\end{abstract}



\end{frontmatter}
%%%%%%%%%%%%%%%%%%%%%%%%%%%%%%%%%%%%%%%%%%%%%%
%% Please use \tableofcontents for articles %%
%% with 50 pages and more                   %%
%%%%%%%%%%%%%%%%%%%%%%%%%%%%%%%%%%%%%%%%%%%%%%
%\tableofcontents

%%%%%%%%%%%%%%%%%%%%%%%%%%%%%%%%%%%%%%%%%%%%%%
%%%% Main text entry area:

\section{Introduction}\label{section_1_intro}

Topological data analysis (TDA) uses multiscale measurements of topology to summarize the \emph{shape} of data. Persistent homology is one of the principal tools in this field, producing a finite descriptor called a persistence diagram that encodes topological and geometric information across scales \cite{comptopEd}. Persistence diagrams have found widespread application in diverse domains including materials science, neuroscience, and computational biology, where shape information provides discriminative features for downstream statistical tasks \cite{vecSurvey}.

In this paper, we study regression where the predictor variable $X$ takes values in the space of persistence diagrams and the response variable $Y$ is real-valued. Given independent observations $(X_1, Y_1), \ldots, (X_n, Y_n)$ from a joint distribution, we wish to estimate the regression function $m(x) = \mathbb{E}[Y \mid X = x]$. This problem arises naturally when one wishes to predict a real-valued quantity such as a physical parameter, a clinical outcome, or a geometric property from topological summaries of complex data.

A fundamental challenge in regression using persistence diagrams is that the space of persistence diagrams, equipped with natural metrics such as the bottleneck or Wasserstein distance, is not a Hilbert space. This precludes direct application of standard regression techniques that rely on vector space structure. The dominant strategy in the literature has therefore been to embed diagrams into a Hilbert space, either explicitly through finite-dimensional vectorization or implicitly through kernel methods \cite{vecSurvey}. 

To this end, several vectorization schemes have been studied. Persistence landscapes \cite{bubenik2015statistical} embed diagrams into a Banach space of functions that can be discretized and Persistence images \cite{adams2017persistence} discretize a weighted Gaussian representation into fixed-length vectors. Kernel methods provide an alternative route via implicit embedding into a reproducing kernel Hilbert space, with kernels such as the persistence scale-space kernel \cite{reininghaus2015stable} and the sliced Wasserstein kernel \cite{carriere2017sliced} enabling kernel regression and Gaussian process regression \cite{kwitt2015statistical}. A variety of other vectorization methods exist, surveyed in \cite{vecSurvey}, though their application to regression problems (such as in \cite{PH_linearModels} or \cite{templates}) has been less systematic. 

Notably, all of these approaches operate by embedding persistence diagrams into a Hilbert space. However, such embeddings cannot preserve the native Wasserstein or bottleneck geometry without distortion \cite{metric_distortion}. Consequently, neighborhood relationships in diagram space may not be faithfully represented in the embedded space, which is particularly relevant for statistical procedures that rely on local metric structure.

%However, authors in \cite{metric_distortion} show that no bi-Lipschitz embedding of persistence diagrams equipped with standard diagram distances into a finite-dimensional Hilbert space exists, even when restricting to diagrams with uniformly bounded cardinality. Moreover, for several stable embeddings into infinite-dimensional Hilbert spaces proposed in the literature, any lower Lipschitz bound must depend explicitly on the cardinality of the diagrams. This highlights the fundamental limitation of vectorizing persistence diagrams, that of metric distortion introduced by the embedding which can obscure the true geometric structure of diagram space, especially in low-sample regimes where such distortions are not mitigated by large amounts of data. 

Furthermore, there are compelling reasons to operate directly within diagram space. The bottleneck and Wasserstein distances underpin much of the statistical theory for persistence diagrams. Foundational stability results \cite{cohen2007stability, cohen2010lipschitz} established that small perturbations in input data induce small changes in persistence diagrams. Mileyko, Mukherjee, and Harer \cite{Mileyko_2011} proved that the space of persistence diagrams equipped with the $p$-Wasserstein metric is a complete, separable metric space, enabling rigorous probability theory on diagrams. This has supported the development of Fr\'{e}chet means and laws of large numbers \cite{KT_Frechet}, confidence sets via bootstrap methods \cite{fasy2014confidence}, null hypothesis testing for comparing distributions of diagrams \cite{robinson2017hypothesis}, and distance covariance methods for measuring dependence between persistence diagrams and response variables \cite{turner2019same}. Hence, various statistical tools have already been developed which operate with respect to the Wasserstein or bottleneck metrics. 

Despite this, to our knowledge, methods for regression from persistence diagrams to real values that operate directly in the native metric space without any Hilbert embedding have been essentially absent from the literature. We propose to address this gap by developing a regression framework directly on the space of persistence diagrams. Since there is no natural parametric model for probability distributions on this space, we adopt a non-parametric local averaging approach, which relies only on the underlying metric structure of the space of predictors. 

In particular, we study kernel-based local averaging estimators, originally proposed by Nadaraya \cite{Nadaraya} and Watson \cite{watson}, which constructs the estimator as an average of the known response values weighted by the proximity of the new input $x$ to the known predictors: 
\[
\hat{m}_n(x) = \frac{\sum_{i=1}^n K_h(d(x, X_i)) Y_i}{\sum_{i=1}^n K_h(d(x, X_i))},
\]
where $K_h$ is a smoothing kernel function with bandwidth $h > 0$ and $d$ is a metric on the predictor space. In the Euclidean case, the Nadaraya–Watson estimator has been extensively studied and shown to be consistent in mean squared error \cite{DW_consistency, SandS_consistency} under suitable conditions. Subsequent work extended these results to infinite-dimensional settings and more generally to regression with predictors taking values in abstract metric spaces. In these broader contexts, mean square consistency has been established under different regularity assumptions placed on the regression function or the probability distributions or the metric space itself \cite{Forzani, Abraham}. Hence, our goal of establishing consistency for our estimator on diagram space can directly leverage this theory.

Another question remains regarding the rate of convergence for such an estimator on our space. In the Euclidean case, if we assume the regression function $m$ to be in the class of functions $\mathcal{F}_p$ where:
\begin{align*}
    \mathcal{F}_{p} = \Big\{ m:[0,1]^d \to \mathbb{R} \ \Big| \ m \text{ is $p$-Hölder smooth} \Big\},
\end{align*}
the optimal minimax rate of convergence is shown to be $O(n^{-2p/(2p+d)})$. That is:
\begin{align*}
    \mathcal{R}_n(\mathcal{F}_{p}) = \inf_{\hat{m}_n} \sup_{m \in \mathcal{F}_{p}} \mathbb{E}\big[\, (\hat{m}_n(X) - m(X))^2 \,\big] \;\; \asymp \;\; O(n^{-2p/(2p+d)})
\end{align*}
where the infimum is taken over all estimators $\hat{m}_n$ (i.e all measurable functions of the data). When $p = 1$, the Hölder class $\mathcal{F}_{1}$ reduces to the class of Lipschitz continuous functions. In this case, the optimal minimax rate of convergence is
\begin{align*}
    \mathcal{R}_n(\mathcal{F}_{1}) \asymp n^{-2/(2+d)}    
\end{align*}
which kernel estimators have been proven to achieve \cite{Gyofri_book, stone1982optimal}. Notably, in the Euclidean setting the convergence rate depends on the ambient dimension $d$ in the denominator, codifying the curse of dimensionality. When predictors take values in an infinite-dimensional or general metric space, the effective “dimension” governing convergence is no longer the ambient Euclidean dimension, but an intrinsic notion of the space’s size or complexity. For instance, rates can depend on small ball probabilities or the Assouad dimension of the predictor space \cite{FerratyBook, knn_F_rate, Kulkarni, daboniang2003, kpotufe2009fast}. In this work, we similarly leverage covering number bounds on a restricted diagram space to establish a corresponding upper bound on the rate of convergence of our kernel estimator.

\subsection{Our Contributions}
Given the discussions above, our main contributions are as follows:
\begin{enumerate}
    \renewcommand{\labelenumi}{\emph{(\roman{enumi})}} % Changes only this list
    \item \emph{Restricted Diagram Space.} First, we introduce a restricted space of persistence diagrams $\mathcal{D}^{N,R}$ corresponding to a finite paradigm that reflects the bounded nature of most practical applications. We establish that this space is compact under the $p$-Wasserstein distance and derive explicit bounds on its covering numbers.
    \item \emph{Asymptotic Results.} Next, we prove consistency of the kernel regression estimator on this metric space and derive upper bounds on the rate of convergence in mean squared error. The rate we obtain is $O(n^{-2/(2N+2)})$, where $N$ bounds the number of points in the persistence diagrams. This reveals that the parameter $2N$ (corresponding to $N$ points with birth-death coordinates) plays the role of dimension in determining the convergence rate, providing a satisfying parallel with classical results where rates depend on the ambient dimension.
    \item \emph{Empirical Validation.} Through computational experiments, we demonstrate strong performance of the proposed estimator. We find that in low-sample regimes where training data is limited, the performance benefit of operating directly in diagram space is most apparent.
\end{enumerate}

\subsection{Outline}
In Section \ref{Background}, we provide background on persistent homology and persistence diagrams. In Section \ref{Diagram_Space_Section} we define the restricted diagram space $\mathcal{D}^{N,R}$ for which we prove compactness, and establish covering number bounds. Then in Section \ref{asymptotic_results_section}, we present the general theory for kernel regression on metric spaces and apply it to obtain consistency and convergence rates for persistence diagram predictors. Finally, in Section \ref{performance_section}, we present experimental results benchmarking our method against existing approaches.

%The second approach attempts to work directly on the space of persistence diagrams, leveraging its metric and geometric structure as in \cite{KT_Frechet}, \cite{aprox_KNN_PD} and \cite{k-means_PD}. Our method takes after the second approach, as we work directly on the metric space of persistence diagrams. 

%\subsection{Outline}
%In Section \ref{Diagram_Space_Section}, we will define a restricted space of all persistence diagrams $\mathcal{D}^{N,R}$ along with a distance metric $d_p$ under which the space $(\mathcal{D}^{N,R}, d_p)$ becomes a metric space. Here the $N$ bounds the number of points in the persistence diagrams and $R$ bounds where these points can lie. We prove this restricted space is compact and provide covering number bounds.

%In Section \ref{asymptotic_results_section} we consider the kernel regression estimators in Definition \ref{KernelEstimatorDefinition} on this metric space. This requires a generalisation of the existing theory where the predictor is a vector. The proofs ar analogous and we include them in an appendix for completeness. In Section \ref{sec:convergence persistence}, we apply this theory to the finite paradigm setting of persistence diagrams as predictor. We establish this estimator as consistent and derive an upper bound on the rate of convergence in MSE. Finally in Section \ref{performance_section}, we benchmark the performance of our estimator against other popular machine learning and TDA methods.

\section{Persistent Homology and Persistence Diagrams}
\label{Background}
Here we will provide a very brief introduction to persistent homology. A more complete coverage can be found in \cite{comptopEd}. Although homology can be computed over any ring, to compute persistent homology we need a field. This is usually $\mathbb{F}_2$ for computational purposes. 

To define persistent homology, we start with a nested sequence of topological spaces,
\begin{align}\label{filteredspace}
A_0 \subseteq A_1 \subseteq A_2 \subseteq \ldots A_n =A.
\end{align}
Often this sequence arises from the sublevel sets of a function, $f \colon A \to \R$, with $A_i =
f^{-1}(-\infty,t_i]$ for a sequence $-\infty< t_0 \leq t_1 \leq \ldots \leq t_n\leq \infty$. 
The sequence induces linear maps on homology for any dimension r:
\[H_r(A_0) \to H_r(A_1) \to \ldots \to H_r(A_n).\]

We are interested in when homology classes appear and disappear in this sequence.
Let $\phi^j_i \colon H_r(A_i) \to H_r(A_j)$ be the linear map on homology induced by the inclusion $A_i\to A_j$. Observe that if $i<j<k$ then $\phi_j^k \circ \phi_i^j = \phi_i^k$. The homology
class $\gamma \in H_r(A_i)$ is said to be born at $A_i$ if it is in the cokernel of $\phi_{i-1}^i$. This same class is said to die
at $X_j$ if its image in $H_r(A_{j-1})$ is not in the image of $\phi_{i-1}^{j-1}$, but its image in $H_r(A_j)$ is in the image of $\phi_{i-1}^j$. In the case that the spaces arose from the level sets of a function $f$ as defined above, we define say that $\gamma$ is born at time $t_i$, dies at time $t_j$, and its lifetime is $t_j-t_i$.

A discrete description of the persistent homology of \eqref{filteredspace} is the multiset of points in the extended plane where we include $(t_i,t_j)$ with multiplicity the dimension of the set of homology classes born at $t_i$ and dying at $t_j$. This is the information recorded as the ($r$th dimensional) persistence diagram. 

Although we have here restricted our definition of persistent homology to that of a finite nested sequence, it can be defined more generally for nested parameterized families of topological spaces $\{A_t\colon t\in \R\}$ with the condition that $A_s\subseteq A_t$ whenever $s\leq t$ (alongside some technical finiteness conditions). We can think of this as a filtered space with filtration parameter $t$. For ``nice'' filtered spaces the persistence diagram can be defined analogously to the finite case. For details see \cite{crawley2015decomposition}.

It is worth observing that although we are computing homology we are often capturing geometric information through the filtration parameter. For example, if we are considering a filtration of $\R^2$ by the distance from a circle of radius $r$ then the $1$st dimensional persistent homology will have exactly one class born at $0$ (the time the circle first appears) and dying at $r$ (the time when the circle is filled in). 

\subsubsection{Motivating Example}
\begin{figure}[H]
  \centering
  \includegraphics[width=1\textwidth]{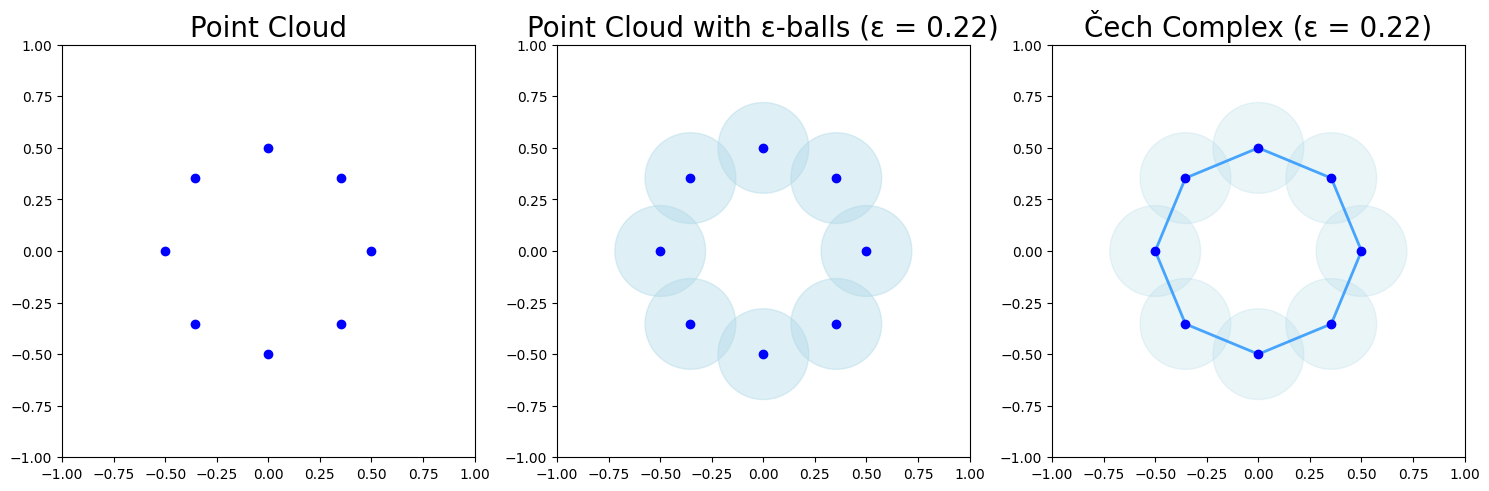}
  \caption{The point cloud $P$ (left), the $\epsilon$-blurred point cloud $\bigcup_{p \in P} B(p, \epsilon)$ (center) and the resulting simplicial complex (right) $\check{C}_\varepsilon(P)$ for $\epsilon = 0.22$ .}
  \label{fig:toy1}
\end{figure}

Consider a point cloud $P \subset [0,1] \times [0,1] \subset \mathbb{R}^2$ in the square portion of the plane as visualized in Figure \ref{fig:toy1}. Our goal is to capture the overall \emph{shape} information in this point cloud which we know to be roughly circular. To do this, we effectively \emph{blur our vision} by placing an $\varepsilon$ ball at each point in $P$, giving us the set $\bigcup_{p \in P} B(p, \varepsilon)$. The relevant information from the blurred point cloud is the overlapping of the $\varepsilon$-balls which connects or blurs the points together in the shape of a loop. This overlap information can be stored in a simplicial complex constructed on $\bigcup_{p \in P} B(p, \varepsilon)$. The C\v ech complex is an example of such a construction: 
\begin{definition}[C\v ech complex]
Let $P = \{p_1, p_2, \dots, p_n\} \subset \mathbb{R}^d$ be a finite point cloud and $\varepsilon > 0$ a scale parameter. The \emph{Čech complex} at scale $\varepsilon$, denoted $\check{C}_\varepsilon(P)$, is the abstract simplicial complex
\[
\check{C}_\varepsilon(P) = \Big\{ \sigma \subseteq X \;\Big|\; \bigcap_{p_i \in \sigma} B(p_i, \varepsilon) \neq \varnothing \Big\}
\]
That is, a subset of points forms a simplex if and only if the corresponding $\varepsilon$-balls have a non-empty intersection.
\end{definition}
Figure \ref{fig:toy1} illustrates the process of constructing the C\v ech complex with $\epsilon = 0.22$, in which we can see a connected loop shape emerge. 

Finally, we can leverage the tools from algebraic topology, namely calculating the simplicial homology of $\check{C}_\varepsilon(P)$. For $\epsilon = 0.22$, we get a single non-trivial $1$-cycle which corresponds to a single generator of the first dimensional homology group $H_1$. On the other hand, if our point cloud dataset was linear, the same method would not yield any $1$-cycles in the corresponding C\v ech complex for any value of $\varepsilon$. Hence, counting the number of generators for $H_1$ in $\check{C}_\varepsilon(P)$ quantifies the circular shape information from $P$ if the value of $\varepsilon$ is well chosen. This allows us to distinguish circular and linearly shaped point clouds.  

\begin{figure}[H]
  \centering
  \includegraphics[width=1\textwidth]{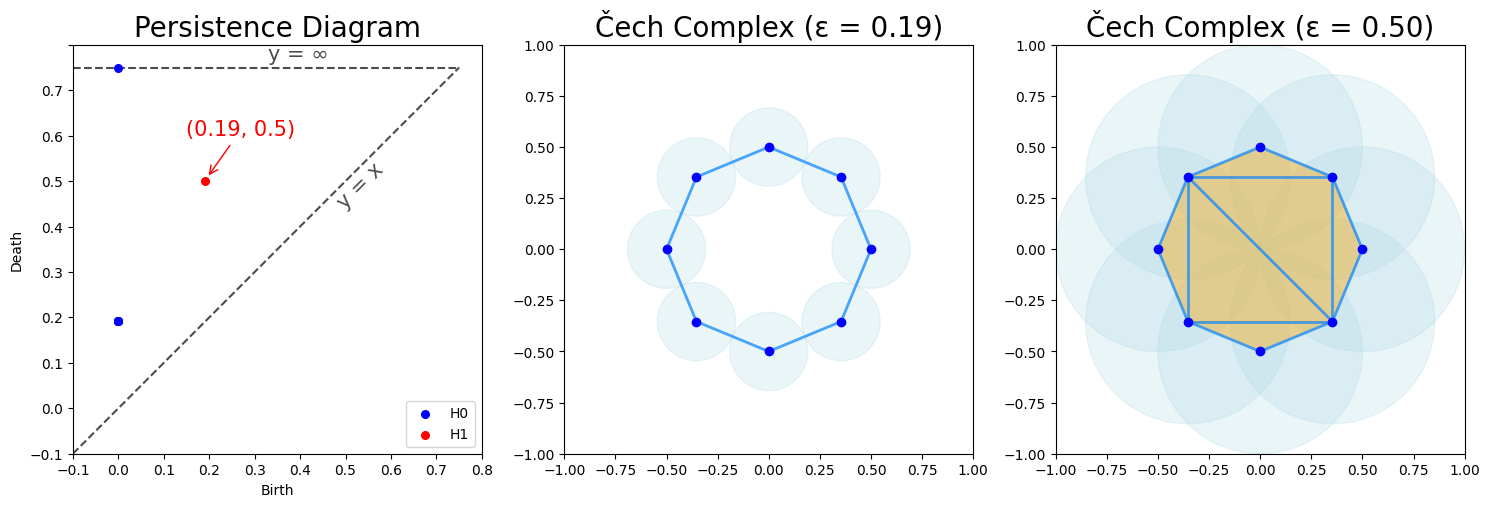}
  \caption{The persistence diagram (left) displays the birth-death co-ordinates for every homological feature as $\epsilon$ varies. A dimension $1$ homology feature appears for $\epsilon = 0.19$ as a loop in the C\v ech complex (center) and disappears at $\epsilon = 0.5$ when the loop is filled in with high order simplices (right)}
  \label{fig:toy2}
\end{figure}

It is apparent that the correct choice of $\varepsilon$ is central for the success of the above method. The key idea behind persistent homology is that instead of choosing a particular value for $\varepsilon$, we calculate the homology information at every value of $\varepsilon$. The persistence diagram stores exactly this information in a slightly different format. As the filtration parameter $\varepsilon$ varies, it record the values of $\varepsilon$ when a new homology feature is born (birth coordinate) and when that feature disappears from the simplicial complex (death coordinate). Hence, a persistence diagram can just be considered as a multiset of points (since multiple homology features can have the same birth-death coordinates) above the diagonal on the plane (since every homology feature must appear before disappearing). Figure \ref{fig:toy2} illustrates the persistence diagram for our point cloud $P$ and how each point reflects the birth and death of a homology feature. 

By varying $\epsilon$ and capturing the changes in homology features, a persistence diagram captures shape information at all scales or all levels of blurring in a single object. Beyond this example, persistence diagrams can be constructed on a variety of complex data types such as images \cite{PH_linearModels} or graphs \cite{graphFiltrations}.

While we have defined the construction of a persistence diagram from a finite nested sequence of topological spaces, it will be important to define the space of persistence diagrams independent of the process of persistent homology.

Let $\R^{2+}=\{(a,b)\in \R^2: a\leq b\}$. This is the part of the plane above the diagonal and points lying here correspond to persistent homology classes with finite lifetimes. Since we wish to also want to consider persistent homology classes of infinite duration we will also want to include lines at infinity. Let $\L_{\infty}:=\{(a,\infty)\colon a\in \R\}$. 
%Persistence diagrams will be multiset of points in $\L_\infty \cup  \R^{2+}$. 
For tractability we will impose some finiteness conditions, namely that only finitely many classes have infinite lifetimes and that the sum of all the finite lifetimes is finite. This restriction is not onerous in applications where generally we have finite sized data as input.

%Any persistence diagram is a multiset of points above the $y = x$ diagonal. Hence, the space of persistence diagrams can be defined as below:
%\begin{definition}[Space of Persistence Diagrams]
%Let $\mathcal{D}$ be the space of persistence diagrams with each diagram $d \in \mathcal{D}$ defined as $d = (S, \Delta, \mu)$ where:
%\begin{itemize}
 %       \item[$\bullet$] $\Delta = \set{ (b,d) \in \mathbb{R}^2: b = d}$ is the diagonal
  %      \item $S \subset \set{ (b,d) \in \mathbb{R}^2: b < d}$ is a countable set of points above the diagonal
   %     \item[$\bullet$] $\mu: S \cup \Delta \to \overline{\mathbb{N}}$ gives the multiplicity of each point in the diagram where:
    %    \begin{align*}
     %       \mu(p) = \begin{cases} 
      %              < \infty & \text{if } p \in S \\
       %             \infty & \text{if } p \in \Delta
%                \end{cases}
 %       \end{align*}
%\end{itemize}
%\end{definition}

\begin{definition}
A \emph{persistence diagram} $X$ is a countable multi-set of $\L_\infty \cup \R^{2+}$ such that 
 the number of elements in $X|_{\L_{\infty}}$ is finite, and
$\sum_{(x_i,y_i)\in X|_{\R^2+}}(y_i-x_i)<\infty$.  
\end{definition}

\begin{remark}
For full generality, we could extend this definition to also include points with $-\infty$ as the birth time. This is rarely useful in applications and will not be needed in our examples later. We removed this extra generality for clearer exposition, but there is no theoretical impediment to include these extra kinds of points. %For each point with $-\infty$ as a birth time we move to a new connected component in the sapce of persistence diagrams
\end{remark}

Let $\D$ denote the space of all persistence diagrams. We will consider a family of metrics which are analogous to the $p$-Wasserstein distances on the space of probability measures and to the $L^p$ distances on the space of functions on a discrete set. $\R^{2+}$ inherits natural $L^p$ distances from $\R^2$. For $p\in [1,\infty)$ we have 
$\|(a_1, b_1) - (a_2, b_2)\|_p^p = |a_1-a_2|^p + |b_1-b_2|^p$ and 
$\|(a_1, b_1) -(a_2, b_2)\|_\infty = \max \{ |a_1-a_2|, |b_1, b_2|\}$.

 Let $\Diag$ denote an abstract element representing the diagonal $\{(x,x)\colon x\in \R\}$. With a slight abuse of notation we write $\|(a, b) - \Delta\|_p $ to denote the shortest $L^p$ distance from $(a,b)$ in to a point in the diagonal set in $\R^2$.  Thus 
\[\|(a, b) - \Delta\|_p = \inf_{t\in \R}\|(a,b)-(t,t)\|_p= 2^{\frac{1}{p}-1}|b-a|\] for $p<\infty$, and $\|(a, b)-\Delta\|_\infty=\inf_{t\in \R}\|(a,b)-(t,t)\|_\infty = |y-x|/2$. $\L_\infty$ also inherits natural $L^p$ distances from the $L^p$ metric on $\R$; $\|(a_1,\infty)-(a,\infty)\|_p=|a_1-a_2|$. We should think of $\L_{-\infty}$ and $\R^{2+}\cup \Delta$ as disjoint parts of a larger space; the distance between any point in $\L_{\infty}$ to any point in $\R^{2+}\cup \Delta$ is infinity.

Distances between persistence diagrams are phrased as optimal transport distances where the allowable transportation plans are called matchings.

\begin{definition}
    Given  persistence diagrams $X$ and $Y$ a \emph{matching} is a triple $(X', Y', \phi)$ with $X'\subset X$, $Y'\subset Y$ and $M:X' \to Y'$ a bijection.
\end{definition}

We define the cost of a matching as a sum of costs of moving points from $X'$ to $Y'$ along $\phi$ and moving the unmatched points to the diagonal. The cost will depend on the choice of $p\in [1,\infty]$.

\begin{definition}
    Given  a matching $(X', Y', \phi)$ between persistence diagrams $X$ and $Y$, its \emph{p-cost}, $c_p(X',Y', M)$ is defined by
    $$c_p(X', Y',\phi)^p=\sum_{x\in X'}\|x-\phi(x)\|_p^p +\sum_{x\in X\backslash X' \cup Y\backslash Y'} \|x-\Delta\|_p^p$$
    for $p\in [1, \infty)$, and 
    $$c_\infty(X', Y',\phi)=\max\{\sup_{x\in X'}\|x-\phi(x)\|_\infty,  \sup_{x\in X\backslash X' \cup Y\backslash Y'} \|x-\Delta\|_\infty\}.$$
\end{definition}

For each $p$ define $$d_p(X,Y) = \inf_{\text{matchings }(X', Y', \phi)} c_p(X', Y', \phi).$$ These distances may be infinite - for example if $X$ and $Y$ contain a different number of points in $\mathcal{L}_\infty$ then $d_p(X,Y)=\infty$ for all $p$.

\section{The Finite Paradigm}
\label{Diagram_Space_Section}
The naked space of persistence diagrams is too large to admit useful properties such as compactness, completeness or separability. Various authors as in \cite{Mileyko_2011}, \cite{templates}, \cite{Bubenik_2018},  \cite{KT_Frechet} impose different restrictions to establish different properties. In this Section, we will define the restricted metric space of persistence diagrams $(\mathcal{D}^{N,R}, d_p)$ which represents the finite paradigm inherent in the finite nature of most input data in applications. We derive theoretical properties of this space required for establishing the asymptotic results in Section \ref{asymptotic_results_section}.

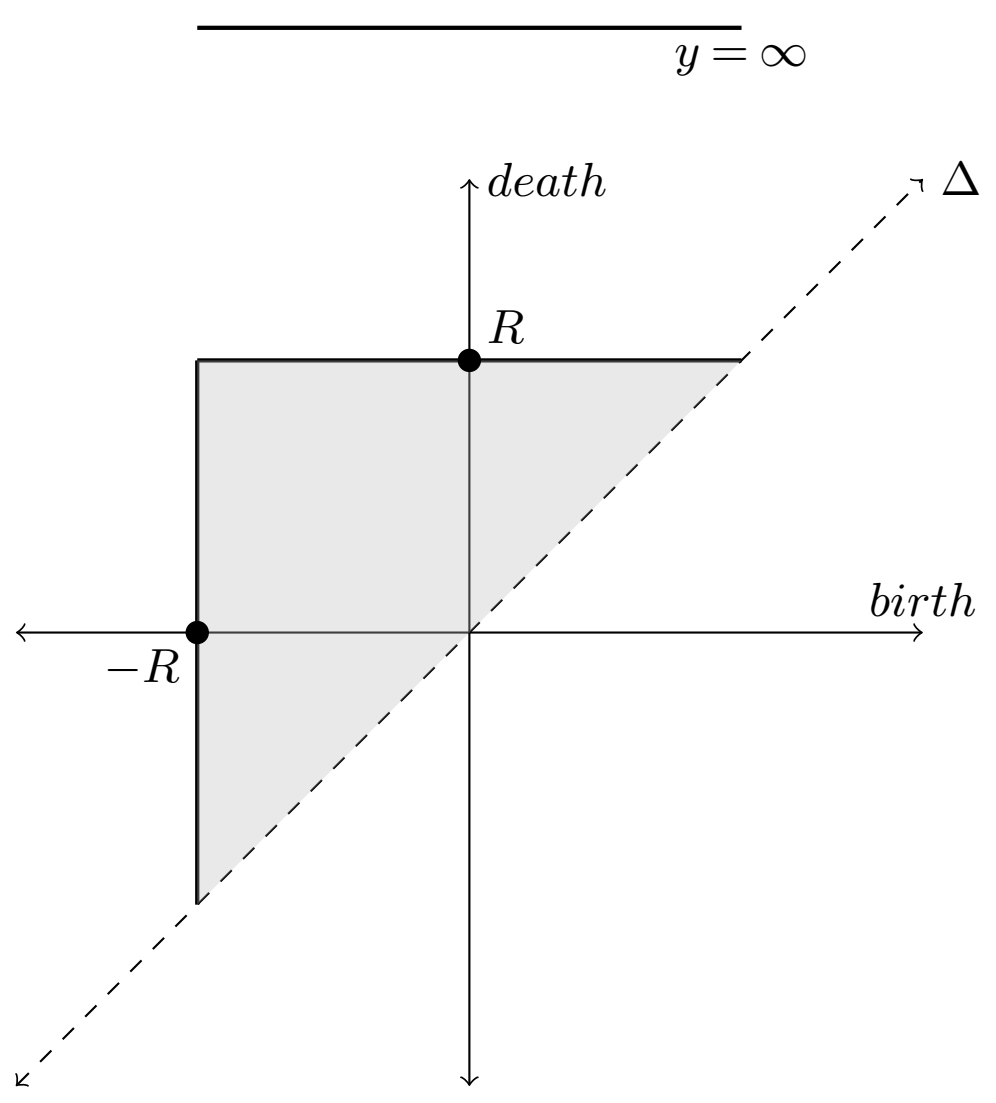
\begin{figure}[H]
      \centering
      \begin{tikzpicture}
            % Axes
            \draw[<->] (-3, 0) -- (3, 0) node[above] {$birth$};
            \draw[<->] (0, -3) -- (0, 3) node[right] {$death$};
            
            % Diagonal line (Delta)
            \draw[<->, dashed] (-3, -3) --(3,3) node[pos=1, right] {$\Delta$};
        
            % Draw region bounds and shade
            \draw[-, thick] (-1.8, 1.8) -- (1.8, 1.8);
            \draw[-, thick] (-1.8, 1.8) -- (-1.8, -1.8);
        
            \fill[gray!50, opacity=0.33] (-1.8, -1.8) -- (-1.8, 1.8) -- (1.8, 1.8) -- cycle;
        
            % Point at (0, 1.8) labeled R
            \draw[fill=black] (0, 1.8) circle (2pt) node[above right] {$R$};
        
            % Point at (-1.8, 0) labeled -R
            \draw[fill=black] (-1.8, 0) circle (2pt) node[below left] {$-R$};
        
            % Line at y = infinity from -R to R
            \draw[-, thick] (-1.8, 4) -- (1.8, 4) node[left, below] {$y = \infty$};
        
          %  \node at (-1, 1) {$T_R$};
        
        \end{tikzpicture}
    
\caption{The set $T_R$ where the points in the persistence diagrams in $\D^{N,R}$ are allowed to be in.}\label{fig:T_R}
\end{figure}

 For our purposes, we will restrict diagram space $\mathcal{D}$ to the finite paradigm as defined below:

\begin{definition}[Finite Paradigm Restricted Space]
    For some $N, R > 0$, let $\mathcal{D}^{N,R}$ be the subset of $\mathcal{D}$ such that for all diagrams $X$ in $\mathcal{D}^{N,R}$ we have:
    \begin{enumerate}[\text{Condition} 1:]
        \item (Finite Total Multiplicity)
         $|X|\leq N$
        \item (Finite Region Bound)
            $x\in X$ implies $x\in [-R, R]\times\{ [-R,R]\cup \{\infty\}\}$ 
    \end{enumerate}
This creates a triangular region and an interval within $\L_\infty$ (depicted in Figure \ref{fig:T_R}) which together we denote $T_R$ which contains all the finitely many points in any persistence diagram.
\end{definition}

By restricting to $\mathcal{D}^{N,R}$
%, we assume that every diagram contains at most $N$ homology features with a global bound $R$ on the birth and death times of those features. 
we get a well-behaved space in which we can avoid many pathologies such as diagrams with infinitely many points accumulating away from the diagonal. We find this restriction is not limiting because persistence diagrams, in most practical settings, arise from studying processes with real world domain bounds which naturally translate into the global bounds on multiplicity and birth-death times. %\textcolor{red}{One feature always has infinite death. How does this fit into the finite paradigm?}\\

For instance, consider a dataset consisting of point clouds in $\mathbb{R}^2$ contained within a bounded square region. Hence, the filtration parameter is necessarily bounded above by the diameter $R$ of the square, providing a natural region bound. Moreover, if the process under study can generate only finitely many points within the square, then the number of homological features across all dimensions is finite, yielding some $N < \infty$ that bounds the number of points in the persistence diagrams. See Section \ref{interpolation_experiment_section} for a concrete experimental example. Similar reasoning applies to persistence diagrams constructed from images (finite image size and resolution) and from graphs (finite number of vertices and edges). In such cases, the finite paradigm is a natural reflection of the underlying data constraints rather than an artificial restriction. 

If persistence diagrams arise from sublevel set filtrations (or from Rips or Cech filtrations of points clouds), then the points in the persistence diagrams with infinite death (those in $\L_\infty$) correspond to the homology classes of the domain and is generally consistent over the dataset of persistence diagrams. This implies we remain in a single connected component within the space of persistence diagrams.

\begin{remark}
  The finite paradigm restriction does have its limitations. For instance, as noted in \cite{templates}, persistence diagrams with infinitely many points may be necessary to encode phenomena such as fractal structures or attractors arising in dynamical systems.
\end{remark}

\noindent Now we establish some properties for $(\mathcal{D}^{N,R}, d_p)$, most importantly compactness and a bound on the covering number. This bound is directly required for bounding the rate of convergence.

\begin{theorem}[Compactness]
\label{compactness_theorem}
For any choice of $p \in [1, \infty]$, the metric space $(\mathcal{D}^{N,R}, d_p)$ is compact.
\end{theorem}
\begin{proof}

Let $(X_n)$ be a sequence of persistence diagrams in $\mathcal{D}^{N,R}$. We can represent each diagram as a vector in $(\mathbb{R}\cup \infty)^{2N}$ since we have at most $N$ points, each with a birth and death coordinate. So we can represent each diagram as the vector $(x_1, ..., x_m, 0, ..., 0)^T$. Note that $0$s used to fill in the vector can be interpreted as the point $(0, 0)$ on the diagonal which has infinite multiplicity. Note that the ordering of points within the persistence diagram in terms of the birth-death coordinate points in the vector representation is arbitrary, but fixed for each diagram. Hence, the same argument will apply for any choice of ordering. \\

Note that the birth coordinate of every point in a diagram $X_n$ is within the range $[-R, R]$ and the death coordinate is in $[-R, R]\cup \{\infty\}$. since $X_n \in \mathcal{D}^{N,R}$. Hence, the vector representation of each diagram is contained within the closed region $([-R, R] \cup {\infty})^{2N}$. Hence the sequence $(X_n)$ can be represented as a sequence of vectors $(x_n)$ in the region $([-R, R] \cup {\infty})^{2N}$. We can look at this sequence of vectors coordinate-wise with the following notation:
\begin{align*}
    (x_n)_{n \in \mathbb{N}} = \big( (x^{(1)}_n, x^{(2)}_n, \dots, x^{(i)}_n, \dots, x^{(2K)}_n)^T \big)_{n \in \mathbb{N}}
\end{align*}
We know $[-R, R]\cup \{\infty\}$  is compact with respect to any $l_p$ metric. Hence, as a finite product of compact sets, we get that $([-R, R] \cup \{\infty\})^{2N}$ is compact with respect to the product topology by Tychonoff's theorem. We also know that the sup-norm metric induces the product topology on $([-R, R] \cup \{\infty\})^{2N}$. Hence, for any choice of $p \in [1, \infty]$, via sequential compactness we get there exists a subsequence $(x_{n_j})$ converging to some $x_* \in ([-R, R] \cup \{\infty\})^{2N}$. That is:
\begin{align}
    \label{proof_subsequence_convergence}
    \max_{1 \leq i \leq 2N} \| x^{(i)}_{n_j} - x^{(i)}_* \|_p  \to 0 \quad \text{as} \quad j \to \infty
\end{align}
We can interpret $x_*$ as a persistence diagram $X_*$ which will be our candidate limit. Note that $X_*$ is in $ \mathcal{D}^{N,R}$ since $x_*$ has $2N$ components where each coordinate is within $[-R, R]\cup \{\infty\}$. Hence, it remains to show $X_{n_j} \to X_*$. For $p \neq \infty$ we have:
\begin{align}
    \label{proof_triangle_inequality}
    d_p(X_{n_j},X_*) %=  \inf_\gamma \left( \sum_{x \in X_{n_j}}\norm{x - \gamma(x)}_{p}^p \right)^{\frac{1}{p}} &
    &\leq \left( \sum_{i = 1}^{2N}\norm{x^{(i)}_{n_j} - x_*^{(i)}}_{p}^p \right)^{\frac{1}{p}}
\end{align}
%where we've chosen the matching $\gamma$ which maps the vector representations of the diagrams as $\gamma(x_{n_j}^{(i)}) = x_*^{(i)}$. 
as $d_p$ is an optimal transport distance and the choice to send $(x_n^{(2j)}, x_n^{(2j+1)})$ to $(x_*^{(2j)}, x_*^{(2j+1)})$ defines a matching.

%Note that the inequality (\ref{proof_triangle_inequality}) holds by the triangle inequality. But then 
By (\ref{proof_subsequence_convergence}) we get:
\begin{align*}
    d_p(X_{n_j},X_*) &\leq \left( \sum_{i = 1}^{2N}\norm{x^{(i)}_{n_j} - x_*^{(i)}}_{p}^p \right)^{\frac{1}{p}}
    \\&\leq (2N)^\frac{1}{p} \left( \max_{1 \leq i \leq 2N} \| x^{(i)}_{n_j} - x^{(i)}_* \|_p \right) \to 0 \quad \text{as} \quad j \to \infty
\end{align*}

Note that $p = \infty$ similarly follows directly from (\ref{proof_subsequence_convergence}). Hence, for any $p \in [1, \infty]$, $(X_n)$ has a converging subsequence, so by sequential compactness $\mathcal{D}^{N,R}$ is compact.
\end{proof}

\begin{definition}[Covering Number]
    Let \( S \) be a subset of the metric space \( (M, d) \). The \( \varepsilon \)-covering number \( \mathcal{C}(S, \varepsilon)\) is the smallest number of open balls of radius \( \varepsilon \) that cover the set \( S \):
    \begin{align*}
        \mathcal{C}(S, \varepsilon) = \inf \left\{ k \geq 1 : \exists \, x_1, \dots, x_k \in S \text{ such that } S \subseteq \bigcup_{i=1}^{k} B(x_i, \varepsilon) \right\}.
    \end{align*}
\end{definition}

\begin{theorem}[Bound on Covering Number]
    \label{CoveringNumberBounds_theorem}
     Let $\epsilon<R/2$. For any choice of the metric $d_p$ for $p \in [1, \infty]$, the covering number of $\mathcal{D}^{N,R}$ as a subset of $(\mathcal{D}, d_p)$ is bounded above by:
    \begin{align*}
       \mathcal{C}(\mathcal{D}^{N,R}, \epsilon) \leq \left({\frac{4R^2N}{\epsilon^2}}\right)^N
    \end{align*}
\end{theorem}
\begin{proof}
As $d_p(X_1, X_2)\leq d_q(X_1,X_2)$ for $p\geq q$ it is sufficient to prove the covering number bound for the case where $p=1$.
To get an upper bound on the covering number $\mathcal{C}(\mathcal{D}^{N,R}, \epsilon)$ in the ambient metric space $(\mathcal{D},d_1)$, we will explicitly construct a finite set $A\subset \mathcal{D}^{N,R}$ such that:
\begin{align*}
    \mathcal{D}^{N,R} \subset \bigcup_{X\in S}B_1(X, \epsilon)
\end{align*}
where $B_1$ is the ball with respect to $d_1$. For this, we must ensure that any diagram $X \in \mathcal{D}^{N,R}$ is within an $\epsilon$ distance of some diagram in $A$.

First we consider the triangular region union a line segment in the extended plane, which we previously denoted $T_R$. We can create a grid of points for the triangular region as depicted in Figure \ref{fig:grid} such that $l_1$ balls with radius $\delta$ cover $T_R$. This grid contains at most 
$$\left(\frac{2\sqrt{2} R}{2\delta}+ 1\right)\left( \frac{\sqrt{2}R}{2\delta}+1\right)$$
locations (with the "$+1$" to account for taking the ceiling). The line segment $[-R, R]\times \{\infty\}$ can be covered by $l_1$ balls with radius $\delta$ regularly spaced $2\delta$ apart - needing at most $\frac{2 R}{2\delta} + 1$. Let $G$ be the set of points in this grid over $T_R$ and the line segment at infinity. %Note that
%$|G|<\frac{2R^2}{\delta^2}$ For sufficiently small $\delta$.
Let $\delta=\frac{\epsilon}{N}$.

\begin{figure}[H]
      \centering
      \begin{tikzpicture}
            % Axes
            \draw[<->] (-3, 0) -- (3, 0) node[above] {$birth$};
            \draw[<->] (0, -3) -- (0, 3) node[right] {$death$};
            
            % Diagonal line (Delta)
            \draw[<->, dashed] (-3, -3) --(3,3) node[pos=1, right] {$\Delta$};
        
            % Draw region bounds and shade
            \draw[-, thick] (-1.85, 1.85) -- (1.85, 1.85);
            \draw[-, thick] (-1.85, 1.85) -- (-1.85, -1.85);
        
            \fill[gray!50, opacity=0.33] (-1.85, -1.85) -- (-1.85, 1.85) -- (1.85, 1.85) -- cycle;
        
            % 45-degree rotated grid of points throughout the triangular region
            % Using rotated coordinates: u (along diagonal), v (perpendicular to diagonal)
            % x = (u - v) / sqrt(2),  y = (u + v) / sqrt(2)
            \foreach \u in {-2.5,-2.2,...,2.8} {
                \foreach \v in {0.0,0.3,...,2.8} {
                    \pgfmathsetmacro{\x}{(\u - \v) / sqrt(2)}
                    \pgfmathsetmacro{\y}{(\u + \v) / sqrt(2)}
                    \pgfmathparse{(\y <= 1.9) && (\x >= -1.9) ? 1 : 0}
                    \ifnum\pgfmathresult=1
                        \draw[fill=blue] (\x, \y) circle (1.5pt);
                    \fi
                }
            }
        
            % Point at (0, 1.8) labeled R
            \draw[fill=black] (0, 1.85) circle (2pt) node[above right] {$R$};
        
            % Point at (-1.8, 0) labeled -R
            \draw[fill=black] (-1.85, 0) circle (2pt) node[below left] {$-R$};
        
            % Line at y = infinity from -R to R
            \draw[-, thick] (-1.85, 4) -- (1.85, 4) node[left, below] {$y = \infty$};
        
            % Dense grid of points on the line at y = infinity
            \foreach \x in {-1.8,-1.5,...,1.8, 1.9} {
                \draw[fill=red] (\x, 4) circle (1.5pt);
            }
        
        \end{tikzpicture}
    \caption{}
    \label{fig:grid}
    \end{figure}

We can show that $T_R$ has a finite cover of the form:
\begin{align*}
    \bigcup_{x_i \in G}{B_{l_1}\left(x_i, \frac{\epsilon}{N}\right)}.
\end{align*}

 Construct the set of persistence diagrams containing at most $N$ points with each of the points lying in $G$, and repetitions are allowed. Denote this set of persistence diagrams $S$.
 % We want to count how many of them there are.  

 Construct the function $\gamma:\R^{2+}\cup \L_{\infty}\to G$ that sends each point to the closest grid point, breaking ties arbitrarily.
% $\set{X_1, \dots, x_N}$ of any $d_i$. It remains however, to account for an arbitrary multiplicity function $\mu$ where $\sum_{x \in d} \mu(x) \leq K$. This is the same as choosing $K$ points from the collection of $N$ points $\set{x_1, ..., x_N}$ with repetition allowed. There are $\dbinom{N + K - 1}{K}$ way to do this, and hence we must construct $L = \dbinom{N + K - 1}{K}$ diagrams where each diagram accounts for a different multiplicity function. Hence, for $p \neq \infty$ and 
For an arbitrary diagram $X \in \mathcal{D}^{N,R}$ consider the diagram where we snap each point to the closest option in $G$. The resulting persistence diagram $\hat{X}=\{\gamma(x)\mid x\in X\}$ is by construction in $A$. By considering the matching which pairs each point to the the closest grid point we have 
\begin{align*}
    d_1(X, \hat{X}) =  \inf_\gamma \left( \sum_{x \in X}\norm{x - \gamma(x)}_1 \right)
    %&\leq 
    %\left( \sum_{i = 1}^N\norm{x - x_i}_{p}^p\right)^\frac{1}{p} \\
    \leq N\frac{\epsilon}{N} = \epsilon
\end{align*}
 This ensures that any point in an arbitrary diagram $X \in \mathcal{D}^{N,R}$ is close to a element in $S$.

%Finally, the smallest $N$ we can pick is $\mathcal{C}(T_R, \frac{\epsilon}{K^\frac{1}{p}})$ by definition of covering numbers. However, this is in the ambient space of the plane ${\mathbb{R}^2}$ with the $l_p$ metric allowing us to determine a more explicit bound as follows:
%\begin{align*}
%    N \leq \mathcal{C}\left(T_R, \frac{\epsilon}{K^\frac{1}{p}}\right) \leq \left( \frac{3K^\frac{1}{p}}{\epsilon}\right)^2 \frac{\text{vol}(T_R)}{\text{vol}(B_{l_p})} \leq \left( \frac{3K}{\epsilon}\right)^2 \frac{\text{vol}(T_R)}{\text{vol}(B_{l_p})}
%\end{align*}
%where $B_{l_p}$ is the unit ball with respect to the $l^p$-norm in $\mathbb{R}^2$. We know $\text{vol}(T_R) = 2R^2$ and vol$(B_p) \leq \text{vol}(B_\infty) = 4$. Hence $N \leq \frac{9R^2K^2}{2\epsilon^2}$. 

Counting the number of elements in $Z$ is the same as the number of ways of choosing at most $N$ points from the collection of $|G|$ points with repetitions allowed. There are $\dbinom{N +|G|}{N}$ way to do this, so 
 $$|Z|= \dbinom{N + |G|}{N}.$$

It remains to show that $|Z|\leq \left({\frac{4R^2N}{\epsilon^2}}\right)^N$. It is trivially true for $N=1$ so we can assume that $N\geq 2$. Recall we have assumed that $\epsilon<R/2$.

First observe that $\dbinom{a+b}{b} \leq \left(\frac{e(a+b)}{b}\right)^b$ for any positive integers $a,b$. 
For $N\geq 2$ and $\epsilon \leq R/2$ we have \begin{align*}
    e(|G|+N) &\leq e\left(\left(\frac{\sqrt{2} RN}{\epsilon}+ 1\right)\left( \frac{RN}{\sqrt{2}\epsilon}+1\right)+ \frac{RN}{\epsilon} +1+N\right)\\
    &\leq e\left(\frac{R^2N^2}{\epsilon^2} + 3.2 \frac{RN}{\epsilon}+N+2\right)\\
    &\leq 4\frac{R^2N^2}{\epsilon^2}.
    \end{align*}
This implies that $|Z|\leq \left(\frac{e(|G|+N)}{N}\right)^N\leq \left(4\frac{R^2N}{\epsilon^2}\right)^N$.

%Combining this inequality and the well known bound for binomial coefficients, we attain our result:
%\begin{align*}
 %   \mathcal{C}(\mathcal{D}^{N,R}, \epsilon) \leq L = \dbinom{|G|+N}}{N}  \\
  %  &\leq \left(\frac{e(|G|+N}{N}\right)^N\\
  %  &\leq \left({\frac{3R^2K}{2\epsilon^2}+1}\right)^K
%\end{align*}
\end{proof}

\noindent Next, if we remove the option of points with infinity as a death time, we can note that our space $\mathcal{D}^{N,R}$ is included in diagram spaces constructed by other authors such as by Mileyko et al. \cite{Mileyko_2011} and Perea et al. \cite{templates}. The conditions imposed by these authors essentially require the number of points that are \emph{far away} from the diagonal to be finite. They define $\mathcal{D}_p$ spaces in analogy with $L^p$ function spaces as below:
\begin{align*}
    \text{Mileyko et al}&: \mathcal{D}_p = \set{X \in \mathcal{D} : d_{p}(X_\varnothing, X) < \infty \text{ for finite } p\in[1,\infty)} \\
    \text{Perea et al}&: \mathcal{D}_\infty = \set{X \in \mathcal{D} : | \, \mathbb{W}^\epsilon(X) \, | < \infty \, \text{ for any } \epsilon > 0}
\end{align*}

    Authors in \cite{templates} have shown the following inclusion of spaces:
    $\mathcal{D}_1 \, \subset \, \mathcal{D}_2 \ ... \ \subset \, \mathcal{D}_{\infty}$
    Since we explicitly have a multiplicity bound $N$, if we excluded points with inifinite death times then would get that $\mathcal{D}^{N,R} \, \subset \, \mathcal{D}_p$ for any $p$.

\begin{remark} 
Authors in \cite{Mileyko_2011} use a slightly different definition for $d_p$, where they fix the $l^\infty$-norm on the $\mathbb{R}^2$ plane instead of letting it vary with $p$. Note this morally doesn't cause issues since for any $p < \infty$ we have 
   $\norm{x - y}_{\infty} \leq \norm{x - y}_{p} \leq 2\norm{x - y}_{\infty}$.
\end{remark}

In \cite{Mileyko_2011}, the space $\mathcal{D}_p$ has been shown to be complete and separable. Again, if we excluded the points with infinite death we could use inclusion of $\mathcal{D}^{N,R}$ as a compact subset in $\mathcal{D}_p$ to show that $\mathcal{D}_p$ is also compact and separable. Since we wish to have up to $N$ points with inifinite death we will instead directly prove completeness and separability.

Finally, we can show that our space is complete and separable.
\begin{theorem}[Completeness and Separability]
    For any choice of $p \in [1, \infty]$, the finite paradigm restricted space $(\mathcal{D}^{N,R}, d_p)$ is complete and separable.
\end{theorem}

\begin{proof}
$(\mathcal{D}^{N,R}, d_p)$  is complete because it is compact ( Theorem \ref{compactness_theorem}). We can prove separability through the bounds on the covering numbers. We can construct a countable dense subset by sequentially including sets of persistence diagrams at the centers of balls that cover at radius $\epsilon_k=1/k$ for $k=\mathbb{N}^+$.
%    Authors in \cite{templates} have shown the following inclusion of spaces:
%    $\mathcal{D}_1 \, \subset \, \mathcal{D}_2 \ ... \ \subset \, \mathcal{D}_{\infty}$
%    Since we explicitly have a multiplicity bound $N$, we get that $\mathcal{D}^{N,R} \, \subset \, \mathcal{D}_p$ for any $p$. In \cite{Mileyko_2011}, the space $\mathcal{D}_p$ has been shown to be complete and separable. The addition of a uniformly bounded number of points with infinite death does not invalidate the arguments used. Via the inclusion of $\mathcal{D}^{N,R}$ as a compact subset in $\mathcal{D}_p$, our space also retains these properties.
\end{proof}

\section{Asymptotic Guarantees}
\label{asymptotic_results_section}
Having set up our metric space of predictors $(\mathcal{D}^{N,R}, d_p)$ over which the local averaging kernel estimator can be defined, we will now establish the mean squared consistency of our estimator, as well as an explicit bound on the rate of convergence. We use the following regression set-up:
\begin{definition}[Regression Setup]
\label{regression_setup_definition}
Suppose $(M,d)$ is a metric space and let $\mathbb{P}$ be a Borel
probability measure on $(M,d)$ with compact support $S^* \subset M$.
Let $\{(X_i,Y_i)\}_{i=1}^n$ and $(X,Y)$ be i.i.d. random elements
taking values in $M \times \mathbb{R}$ with common distribution,
satisfying the regression model $Y=m(X)+\epsilon$, where
\begin{itemize}
    \item $m:M\to\mathbb{R}$, $m(X)=\mathbb{E}(Y\mid X)$,
    is continuous;
    \item $\epsilon$ is random noise such that
    $\mathbb{E}(\epsilon \mid X)=0$ and
    $\mathbb{E}(\epsilon^2 \mid X) = \sigma^2<\infty$.
\end{itemize}
\end{definition}

\begin{remark}
\label{reg_setup_justification}
Note that the assumption of $\mathbb{P}$ having compact support may appear restrictive; however, the ambient space of interest $(\mathcal{D}^{N,R}, d_p)$ is itself compact. Furthermore, continuity of $m$ is a mild strengthening of the Besicovitch condition, which is known to be essentially necessary for the consistency of local averaging methods (see \cite{Forzani, knn_inf_dim, daboniang2003}). Thus, these additional assumptions on the classical regression set up are natural in our context.
\end{remark}

We can construct an estimator $\hat{m}$ for the regression function $m$ via a local averaging approach
\[
\hat{m}_n(X) = \sum_{i=1}^n W_{i}(X) Y_i
\]
where the weight $W_i(X)$ is a measurable function of $(X, X_1, \dots, X_n)$. To prove consistency, we want a generalization of Stone's theorem \cite{Stone} which provides a set of conditions on the weighting functions which guarantee consistency in mean squared error. Instead of working with the already generalized work in \cite{Forzani} or \cite{Abraham}, we will work with the Euclidean proof in \cite{Gyofri_book} which is cleaner and more directly applicable under our assumptions and has not been generalized this way before. As the proof is very analogous, we have provided it in the Appendix \ref{proof_for_stones}.

\begin{theorem}[Generalized Stone's Theorem]\label{thm:Stones} Under the regression set up in Definition \ref{regression_setup_definition}, suppose the weight functions $\set{W_i(X)}_{i=1}^n$ of the local averaging estimator satisfy the following conditions:
\begin{enumerate}

% We can simplify the argument and remove the following condition since we have 
%    \item There is a constant $c$ such that for every non-negative measurable function $f$ satisfying $\mathbb{E}f(X) < \infty$ and for any $n$:
%    \begin{align*}
%        \mathbb{E}\left( \sum_{i=1}^{n} |W_{i}(X)| f(X_i) \right) \leq c \, \mathbb{E}f(X).
%    \end{align*}
    \item There is a constant $D \geq 1$ such that for all $n$:
    \begin{align*}
    \mathbb{P} \left[ \sum_{i=1}^{n} |W_{i}(X)| \leq D \right] = 1
    \end{align*}
    \item For all $\delta > 0$:
    \begin{align*}
    \lim_{n \to \infty} \mathbb{E} \left[ \sum_{i=1}^{n} |W_{i}(X)| \cdot \mathds{1}_{\{d(X_i, X) > \delta \}} \right] = 0
    \end{align*}
    \item 
    \begin{align*}
    \lim_{n \to \infty} \sum_{i=1}^{n} W_{i}(X) = 1 \quad \text{in probability}.
    \end{align*}
    \item 
    \begin{align*}
    \lim_{n \to \infty} \mathbb{E} \left[ \sum_{i=1}^{n} W_{i}(X)^2 \right] = 0
    \end{align*}
\end{enumerate}
Then $\hat{m}_n(X)$ is consistent in MSE, that is:
    \begin{align*}
        \mathbb{E}\left( \, | \hat{m}_n(X) - m(X) |^2 \, \right) \to  0 \quad \text{as} \quad n \to \infty
    \end{align*}
\end{theorem}

The kernel estimator is a type of local averaging estimator and so we can define it in terms of weight functions as below:
\begin{definition}[Regular Kernel Function]
A function $K: \mathbb{R} \to [0,\infty]$ which satisfies
$b \cdot \mathds{1}_{[-r, r]}(x) \leq K(x) \leq \mathds{1}_{[-R, R]}(x)$ for some $0 < r \leq R$ is referred to as a regular kernel function.
\end{definition}

\begin{definition}[Kernel Regression Estimator]
The kernel regression estimator $\hat{m}_n(X) = \sum_{i=1}^n W_i(X, X_1, ..., X_n) Y_i$ is defined with the weighting functions:
\begin{align*}
    W_i(X, X_1, ..., X_n) = \frac{K \big( \frac{d(X_i, X)}{h_n} \big)}{\sum_{j=1}^n K \big( \frac{d(X_j, X)}{h_n} \big)}\mathds{1}_{A}
\end{align*}
where $A = \left\{ \omega \in M^{n+1} : {\sum_{j=1}^n K \big( \frac{d(X_j, X)}{h_n} \big)} \neq 0 \right\}$ for some choice of regular kernel function and choice of bandwidth parameter $h_n$.
\label{KernelEstimatorDefinition}
\end{definition}

It remains to verify that the kernel weights satisfy the conditions of Theorem \ref{thm:Stones}. This proof is a direct generalization of the proof for the Euclidean case in \cite{Gyofri_book} (cf. Theorem 5.1), hence it is relegated to the Appendix \ref{proof_for_consistency}.

\begin{theorem}[Consistency]
    \label{consistency_theorem}
    Consider the kernel estimator in Definition \ref{KernelEstimatorDefinition} under the regression set up in Definition \ref{regression_setup_definition}. Suppose $(h_n)$ is a sequence of positive numbers where $h_n \to 0$ as $n \to \infty$ such that 
        \begin{align*}
            \frac{\mathcal{C}(S^*,rh_n/2)}{n} \to 0 \text{ as } n \to \infty
        \end{align*}
    where $r$ is the inner support radius of the chosen regular kernel. Then we have $\hat{m}_n(X) \to m(X)$ as $n \to \infty$ in the $L^2$-norm. Or equivalently:
    \begin{align*}
        \mathbb{E}\left( \, | \hat{m}_n(X) - m(X) |^2 \, \right) \to  0 \quad \text{as} \quad n \to \infty
    \end{align*}
\end{theorem}

\begin{remark}
\label{hn_existence_remark}
There always exists a sequence $(h_n)$ satisfying the bandwidth condition of Theorem \ref{consistency_theorem}. Since $S^*$ is compact, its covering number $\mathcal C(S^*,\epsilon)$ is finite for every $\epsilon>0$. Hence, one may choose $h_n\to0$ sufficiently slowly so that
\[
\frac{\mathcal C(S^*,rh_n/2)}{n}\to0.
\]
For example, one may choose $(h_n)$ such that
\[
\mathcal C(S^*,rh_n/2)\leq n^{1/2}
\]
for all sufficiently large $n$. Such a sequence exists because the covering number is finite at every positive radius, and the bandwidth can be chosen to decrease arbitrarily slowly.
\end{remark}

Note that Theorem \ref{consistency_theorem} can be directly applied to our metric space of interest $(\D^{N,R}, d_p)$ given Remark \ref{reg_setup_justification}. Hence, it remains to derive a bound on the rate of convergence for our estimator in line with our second objective. To achieve this, we adapt and extend the Euclidean proof from \cite{Gyofri_book}, together with the explicit covering number bound established in Theorem \ref{CoveringNumberBounds_theorem}.
\begin{theorem}[Explicit Bound on Rate of Convergence]
    \label{ROC_bound_theorem_applied}
    Consider the kernel estimator in Definition \ref{KernelEstimatorDefinition} under the regression set up in Definition \ref{regression_setup_definition} where $(M,d) = (\mathcal{D}^{N,R},d_p)$. Suppose:
    \begin{enumerate}
\begin{comment}
    %This is included directly in the regression set up now
    \item $\forall x\in M$, we have $\mathbb{E}\left( \, \epsilon^2 \, | \, X = x\,\right) < \sigma^2_{U} < \infty$    
\end{comment}
    \item For some $C > 0$, $m$ is Lipschitz continuous:
    \begin{align*}
        |m(x) - m(y)| \leq C \cdot d_p(x,y)
    \end{align*}

\begin{comment}
%We can state the bandwidth sequence explicitly for simplification now
    \item $(h_n)$ is a sequence of positive numbers where $h_n \to 0$ as $n \to \infty$ such that 
        \begin{align*}
            \frac{\mathcal{C}(S^*,h_n/2)}{n} \to 0 \text{ as } n \to \infty    
        \end{align*}
\end{comment}
    
    \item The kernel function $K(x)$ is the naive kernel
    \begin{align*}
        K(x) = \mathds{1}_{[-1,1]}(x)
    \end{align*}
    \item The bandwidth sequence is chosen as:
    \begin{align*}
        h_n = n^{-\frac{1}{2N+2}}
    \end{align*}
\end{enumerate}
Then for large enough $n$:
    \begin{align*}
        \mathbb{E}\left[ \, | \hat{m}_n(X) - m(X) |^2 \, \right] \leq a \cdot n^{-\frac{2}{2N + 2}}
    \end{align*}
where $a$ is some constant independent of $n$. Thus, we obtain an upper bound of order $O(n^{-{2}/(2N + 2)})$ on the rate of convergence of the kernel estimator over our restricted diagram space $(\mathcal{D}^{N,R}, d_p)$.
\end{theorem}

\begin{proof}We begin by noting that:
\begin{align*}
    \mathbb{E}\left[ \, ( \hat{m}_n(X) - m(X) )^2 \, \right] &= \int_{S^*}{\mathbb{E} \left[ \, ( \, \hat{m}_n(x) - m(x)\,)^2 \right] } \, d\mathbb{P} \\ &= \int_{S^*}{\mathbb{E} \left[ \mathbb{E}\left[ \, ( \, \hat{m}_n(x) - m(x)\,)^2 \mid X_1, ..., X_n \, \right]\right]} \, d\mathbb{P}
\end{align*}
where the integral is being taken the new input $X$. We will first bound the inner expectation term. Note that every random variable $Z$ with a finite second moment satisfies the identity:
\begin{align*}
\mathbb{E}\big[(Z - a)^2\big] = \mathbb{E}\big[ \, (Z - \mathbb{E}(Z))^2 \, \big] + (\mathbb{E}[Z] - a)^2
\end{align*}
where $a$ is a constant. We can recover the conditional bias and variance decomposition by substituting $Z = \hat{m}_n(x)$ and $a = m(x)$ with the expectation being taken conditionally over the random sample $X_1, ..., X_n$: 
\begin{align*}
    \mathbb{E}\left[ \, | \, \hat{m}_n(x) - m(x)\,|^2 \mid X_1, ..., X_n \, \right] = \mathbb{E}\left[ \, | \, \hat{m}_n(x) - \hat{m}_C(x)\,|^2 \mid X_1, ..., X_n \, \right] + |\, \hat{m}_C(x) - m(x) \,|^2
\end{align*}
where
\begin{align*}
    \hat{m}_C(x) = \mathbb{E} \left[ \, \hat{m}_n(x) \mid X_1, ..., X_n \, \right] = \frac{\sum_{i=1}^n m(X_i)\mathds{1}_{\set{X_i \in B(x,h_n)}}}{\sum_{i=1}^n\mathds{1}_{\set{X_i \in B(x,h_n)}}}\mathds{1}_{A}
\end{align*}

Now we can control each of the conditional bias and variance terms separately. For the first conditional bias term we have:
\begin{align*}
    \mathbb{E}&\left[ \, | \, \hat{m}_n(x) - \hat{m}_C(x)\,|^2 \mid X_1, ..., X_n \, \right] \vphantom{\Bigg{ | }}  \\
    &= \mathbb{E}\left[ \, \left(  \, \frac{\sum_{i=1}^n \left(Y_i - m(X_i)\right)\mathds{1}_{\set{X_i \in B(x,h_n)}}} {\sum_{i=1}^n\mathds{1}_{\set{X_i \in B(x,h_n)}}} \mathds{1}_{A} \, \right)^2 \mid X_1, ..., X_n \, \right] \vphantom{\Bigg{ | }} \\ 
    &=  \frac{\sum_{i=1}^n \mathbb{E}\left( \, \epsilon_i^2 \, | \, X_i \, \right)\mathds{1}_{\set{X_i \in B(x,h_n)}}} {\left( \sum_{i=1}^n\mathds{1}_{\set{X_i \in B(x,h_n)}}\right)^2}\mathds{1}_{A} \vphantom{\Bigg{ | }} \quad \text{since for } a \neq b, \, \, \mathbb{E}(\epsilon_a\epsilon_b \, | \, X_a,X_b) = 0\\
    & \leq \sigma^2 \frac{\sum_{i=1}^n \mathds{1}_{\set{X_i \in B(x,h_n)}}} {\left( \sum_{i=1}^n\mathds{1}_{\set{X_i \in B(x,h_n)}}\right)^2}\mathds{1}_{A} \vphantom{\Bigg{ | }}\\
    &=\sigma^2 \frac{\mathds{1}_{A}} {\sum_{i=1}^n\mathds{1}_{\set{X_i \in B(x,h_n)}}}
\end{align*}
Then the second term is controlled since:
\begin{align*}
    \left(\hat{m}_C(x) - m(x) \right)^2 &= \left( \frac{ \sum_{i=1}^n m(X_i) \mathds{1}_{\set{X_i \in B(x,h_n)}}} {\sum_{i=1}^n\mathds{1}_{\set{X_i \in B(x,h_n)}}}\mathds{1}_A - \left[m(x)\mathds{1}_A + m(x)\mathds{1}_{A^c}\right]\right)^2 \vphantom{\Bigg{ | }}\\
    & \leq \left( \frac{ \sum_{i=1}^n | m(X_i) - m(x) | \, \mathds{1}_{\set{X_i \in B(x,h_n)}}} {\sum_{i=1}^n\mathds{1}_{\set{X_i \in B(x,h_n)}}}\mathds{1}_{A} + |m(x)|\mathds{1}_{A^c} \right)^2 \vphantom{\Bigg{ | }}\\
    & \leq \left( \frac{ \sum_{i=1}^n \left( C h_n \right) \mathds{1}_{\set{X_i \in B(x,h_n)}}} {\sum_{i=1}^n\mathds{1}_{\set{X_i \in B(x,h_n)}}}\mathds{1}_{A} + |m(x)|\mathds{1}_{A^c} \right)^2 \vphantom{\Bigg{ | }}\\
    & \hspace{15em} \text{since } |m(X_i) - m(x)| \leq C \, d_p(X_i, x) \leq C h_n \\
    & \leq C^2h_n^2 \mathds{1}_{A} + m(x)^2 \mathds{1}_{A^c} \text{ where the cross terms above vanish since } \mathds{1}_A \mathds{1}_{A^c} = 0\\
    &\leq C^2h_n^2 + m(x)^2 \mathds{1}_{A^c}
\end{align*}
These two terms can then be used to control the $L^2$-norm as follows:
\begin{align*}
    &\mathbb{E}\left[ \, | \hat{m}_n(X) - m(X) |^2 \, \right]\\
    &= \int_{S^*}{ \mathbb{E}\left[\  \mathbb{E}\left[\ | \, \hat{m}_n(x) - \hat{m}_C(x)\,|^2 \, \mid X_1, ..., X_n \right] \right]} \, d\mathbb{P} \,  +  \int_{S^*}{\mathbb{E}\left[\, |\, \hat{m}_C(x) - m(x) \,|^2 \, \right] } \, d\mathbb{P}  \\
    &\leq \sigma^2 \int_{S^*}  \mathbb{E} \left[ \frac{\mathds{1}_{A}} {\sum_{i=1}^n\mathds{1}_{\set{X_i \in B(x,h_n)}}} \, \right] \, d\mathbb{P} + \int_{S^*} C^2h_n^2 + m(x)^2 \, \mathbb{E}[\mathds{1}_{A^c}] \, d\mathbb{P}
\end{align*}
Noting that $\sum_{i=1}^n\mathds{1}_{\set{X_i \in B(x,h_n)}}$ is a binomial random variable as a sum of i.i.d bernoulli random variables, we can apply Lemma \ref{lemma:binomial_inverse} to upper bound the first integral as:
\begin{align*}
    &\leq \frac{2\sigma^2}{(n+1)} \int_{S^*} \frac{1}{\mathbb{P}(B(x,h_n))} \, d\mathbb{P} + C^2h_n^2 + \int_{S^*} (1 - \mathbb{P}(B(x,h_n)))^n \, m(x)^2 \, d\mathbb{P}
\end{align*}

Since $m$ is continuous on a compact support $S^*$, it attains the extremum value $B = \sup_{x \in S^*}{|m(x)|^2}$. Using this we can control the second integral as follows:
\begin{align*}
    \int_{S^*} (1 - \mathbb{P}(B(x,h_n)))^n \, m(x)^2 \, d\mathbb{P} &\leq B \int_{S^*} (1 - \mathbb{P}(B(x,h_n)))^n  \, d\mathbb{P} \\
    &\leq B \int_{S^*} \text{exp}\Bigl(-n \cdot \mathbb{P}\bigl(B(x,h_n)\bigr) \Bigr) \, d\mathbb{P} \\
    & \hspace{6em} \text{since} \, \, \forall x \in [0,1]: (1-x)^n \leq e^{-nx} \\
    &= B \int_{S^*} \frac{n \cdot \mathbb{P}\bigl(B(x,h_n)\bigr)}{n \cdot \mathbb{P}\bigl(B(x,h_n)\bigr)} \text{exp}\Bigl(-n \cdot \mathbb{P}\bigl(B(x,h_n)\bigr) \Bigr) \, d\mathbb{P} \\
    & \leq \frac{B}{n}\int_{S^*} \frac{1}{ \mathbb{P}\bigl(B(x,h_n)\bigr)} \, d\mathbb{P} \hspace{1em} \text{since} \, \, \forall x \in [0,1]: nxe^{-nx} \leq 1
\end{align*}

\begin{comment}
%No need to decompose the integral now since S* is already compact!

Now let $S = B(d_\varnothing, c)$ be a ball of radius $c$ centered at the empty $d_\varnothing$ (this is just an arbitrary origin). We can use $S$ to decompose the integral as follows:

\begin{align*}
    &= \int_{S \cap S^*} \text{exp}\Bigl(-n \cdot \mathbb{P}\bigl(B(x,h_n)\bigr) \Bigr) \, d\mathbb{P} + \int_{S^c \cap S^*} \text{exp}\Bigl(-n \cdot \mathbb{P}\bigl(B(x,h_n)\bigr) \Bigr) \, d\mathbb{P} \phantom{\Bigg|} \\
    &\leq \int_{S \cap S^*} \text{exp}\Bigl(-n \cdot \mathbb{P}\bigl(B(x,h_n)\bigr) \Bigr) \, d\mathbb{P} + \mathbb{P}\bigl( S^c \cap S^*) \quad \text{since} \, \, \forall x \in [0,1]: e^{-nx} \leq 1
\end{align*}

By making the radius of $S$, $c$ bigger, we can always control $\mathbb{P}\bigl( S^c \cap S^*)$ to be arbitrarily small. The first term can be upper bounded as below:
\begin{align*}
    \begin{split}
    \int_{S \cap S^*} \text{exp}\Bigl(-n \cdot \mathbb{P}\bigl(B(x,h_n)\bigr) \Bigr) \, d\mathbb{P} &= \int_{S \cap S^*} \frac{n \cdot \mathbb{P}\bigl(B(x,h_n)\bigr)}{n \cdot \mathbb{P}\bigl(B(x,h_n)\bigr)} \text{exp}\Bigl(-n \cdot \mathbb{P}\bigl(B(x,h_n)\bigr) \Bigr) \, d\mathbb{P} \\
    & \leq \frac{1}{n}\int_{S \cap S^*} \frac{1}{ \mathbb{P}\bigl(B(x,h_n)\bigr)} \, d\mathbb{P} \hspace{1em} \text{since} \, \, \forall x \in [0,1]: nxe^{-nx} \leq 1
    \end{split}
\end{align*}

\end{comment}
Hence, we get:
\begin{align*}
\mathbb{E}\left[ \, | \hat{m}_n(X) - m(X) |^2 \, \right] \leq \left(\frac{2\sigma^2}{n+1} + \frac{B}{n} \right)\int_{S^*} \frac{1}{\mathbb{P}(B(x,h_n))} \, d\mathbb{P} + C^2h_n^2
\end{align*}
Now applying the small ball probability control Lemma \ref{SBP_control_compact} and simplifying we get the following inequality in terms of $h_n$:
\begin{align*}
\mathbb{E}\left[ \, | \hat{m}_n(X) - m(X) |^2 \, \right]\leq \left(\frac{2\sigma^2 + B}{n} \right)\mathcal{C}\left(S^*, \frac{h_n}{2}\right) + C^2h_n^2
\end{align*}

Now we use the explicit upper bound on the covering number function from Theorem \ref{CoveringNumberBounds_theorem} to construct an explicit $h_n$ sequence. Substituting in the bound we get:
\begin{align*}
    \mathbb{E}\left[ \, | \hat{m}_n(X) - m(X) |^2 \, \right] &\leq \left(\frac{2\sigma^2 + \sup_{x \in S^*}{|m(x)|^2}}{n} \right)\left(16\frac{R^2N}{h_n^2}\right)^N + C^2h_n^2 \\
    &= \left( 16R^2N \right)^N \left(2\sigma^2 + \sup_{x \in S^*}{|m(x)|^2}\right)\left({\frac{1}{nh_n^{2N}}}\right) + C^2h_n^2 \\
    &\leq a \cdot \frac{1}{nh_n^{2N}} + C^2h_n^2
\end{align*}
where $a$ is some constant independent of $n$. Now we can choose our sequence of $h_n$ to be given by:
\begin{align*}
    h_n = n^{-\frac{1}{2N + 2}}
\end{align*}
which we have picked so that we get the following simplification after substituting $h_n$ into the upper bound:
\begin{align*}
    \mathbb{E}\left[ \, | \hat{m}_n(X) - m(X) |^2 \, \right] &\leq a \cdot \frac{1}{nh_n^{2N}} + C^2h_n^2 \\
    &= a \cdot n^{-1} \cdot n^{\frac{2N}{2N + 2}} + C^2 \cdot n^{-\frac{2}{2N + 2}} \\
    &= a' \cdot n^{-\frac{2}{2N + 2}} \quad \text{for large } n
\end{align*}
where $a'$ is just some constant independent of $n$.
\end{proof}

Now there is an interesting point of comparison to be made with the Euclidean case where the optimal rate of convergence has upper bounded of order $O(n^{-{2}/(d + 2)})$ \cite{Gyofri_book} for Lipschitz continuous regression functions in ambient dimension $\mathbb{R}^d$. Thus, the ambient dimension parameter $d$ is effectively replaced by $2N$ in our setting where $N$ is the bound on the number of off-diagonal points in the persistence diagrams. We can interpret the value $2N$ as $N$ points living in $2$-dimensional space (with birth-death coordinates). Intuitively, letting each birth and death coordinate of $N$ points be a free variable, we end up in $\mathbb{R}^{2N}$. Hence, in this context, we can see how the bound on the total multiplicity of a set of persistence diagrams fulfills a dimensionality criteria.

\iffalse
\begin{remark}
It is worth noting there is a slight but important technical difference to our estimator constructor than in the set up in \cite{Forzani}. They have an assumption that the sum of weights in their weighted average is always one, whereas we have a sum of weighted averages as zero outside the union of balls centered at the points. While in the limit we have coverage, and so these coincide, this slight variation would cause frustration if we tried use their theory directly. \textcolor{red}{Confused by this remark as our weights always sum to 1??}
\end{remark} 
\fi

\section{Experimental Results}
\label{performance_section}
Finally, we will test the performance of our estimator. 

\subsection{Implementation Outline}
\label{implementation_steps}
Given the sample data $\set{ \,(x_i, y_i) \,}_{i = 1}^N$ and a new input $x_{N+1}$, we have implemented the following algorithm in Python:
\begin{enumerate}[Step 1:]
    \item For each sample input $x_i$ and the new input $x_{N+1}$, we construct the corresponding set of persistence diagram, resulting in a set of persistence diagrams $\set{d_1, ..., d_N, d_{N+1}}$. We have used the \texttt{Ripser} library \cite{Bauer2021Ripser} which constructs a computationally efficient Vietoris-Rips simplicial complex.
    \item Compute the distances between the new input persistence diagram $d_{N+1}$ and the sample diagrams $\{d_1,\ldots,d_N\}$ using the metric $d_p$; where we take $p=2$ throughout. Distances are computed using the \texttt{GUDHI} library \cite{GUDHI} to ensure computational efficiency. Note that essential classes (infinite death time) are retained and matched only to one another by default. The distances for the $H_0$ and $H_1$ persistence diagrams are evaluated separately and then summed, yielding a collection of distances $\{r_1,\ldots,r_N\}$.
    \item Using the distances, we calculate a set of weights for each input sample $\set{w_1, ..., w_N}$ using Definition \ref{KernelEstimatorDefinition}. We've chosen the Gaussian Kernel $K_h\!\left(z\right) = \frac{1}{\sqrt{2\pi}} \exp\!\left(-\frac{z^2}{2h^2}\right)$ and chosen the bandwidth parameter $h$ via leave-one-out cross-validation. Note that we have not chosen a value for $R$ in Definition \ref{KernelEstimatorDefinition} so we never run into the case where the sum of kernel weights is $0$.
    \item Create the regression estimate by the weighted average of target values $\set{y_i}_{i=1}^N$ using the previously calculated weights.
\end{enumerate}

The Python code used to implement the method, along with all datasets used in the experiments have been made publicly available at \href{https://github.com/saksham-aggarwal-2002/kernel-regression-persistence-diagrams.git}{kernel-regression-persistence-diagrams.git}

\subsection{Interpolation Parameter Experiment}
We begin with demonstrating our method on a toy dataset.
\label{interpolation_experiment_section}

\subsubsection{Dataset}
We illustrate our method on a synthetic toy dataset designed to interpolate continuously between two simple geometric structures. Fix $20$ equidistant angles $\theta_1, \dots, \theta_{20} \in [0, 2\pi]$. Let $U_i^{(x)}, U_i^{(y)} \sim \mathrm{Uniform}(-0.2, 0.2)$ be independent noise variables for each point and each coordinate. For a given interpolation parameter \(t \in [0,1]\), we generate a point cloud \(\{(x_i, y_i)\}_{i=1}^{20} \subset \mathbb{R}^2\) according to
\[
x_i = \cos(\theta_i) + U_i^{(x)}, \qquad
y_i = (1 - t)\cos(\theta_i) + t\sin(\theta_i) + U_i^{(y)}.
\]
When $t = 0$, the point cloud concentrates around a noisy circle, while for $t = 1$ it lies approximately along a noisy line. Intermediate values of $t$ yield a smooth geometric interpolation between these two regimes. The regression task is to predict the parameter $t$ from a noisy point cloud realization. Figure ~\ref{fig:toy3} illustrates representative samples for several values of $t$.

\subsubsection{Methods}
\begin{figure}[H]
  \centering
  \includegraphics[width=1\textwidth]{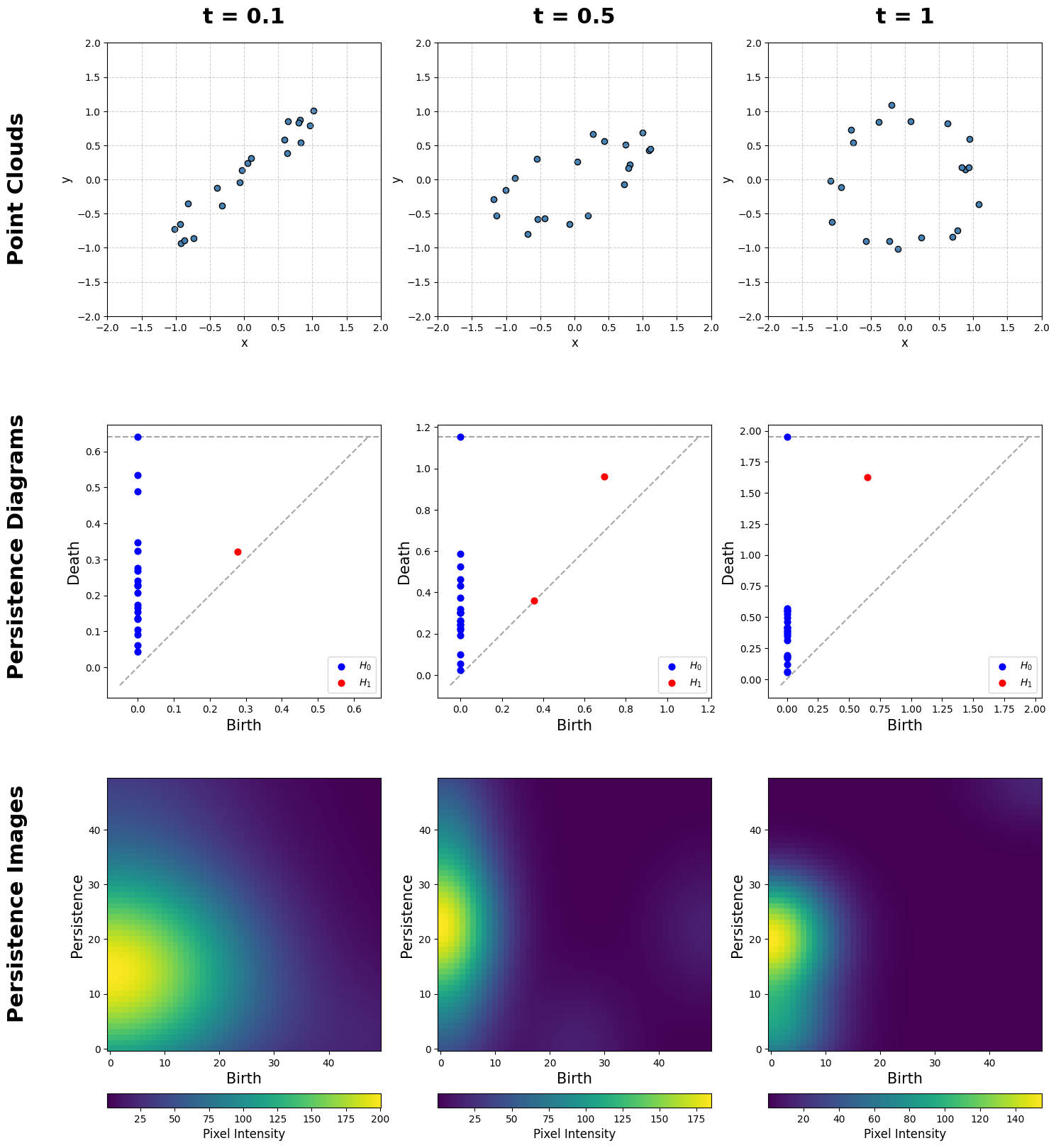}
  \caption{Sample point clouds from the Interpolation Parameter Experiment (\(t = 0\), \(0.5\), and \(1\)) and their associated full persistence diagrams $\mathrm{PD}$ and persistence images $\mathrm{PI}$.}
  \label{fig:toy3}
\end{figure}

To assess the benefits of operating directly in the space of persistence diagrams, we compare our kernel regression method described in Section~\ref{implementation_steps} with a commonly used vectorization approach based on persistence images. In both cases, we employ the same kernel regression framework. The difference lies solely in the representation of the persistence diagrams: either as elements of persistence diagram space equipped with the Wasserstein distance, or as finite-dimensional vectors in Euclidean space via persistence images. Figure~\ref{fig:toy3} includes the persistence diagrams and persistence images associated with the sample point clouds for each value of $t$.

Persistence images provide a stable vectorization of persistence diagrams by convolving each persistence point with a Gaussian kernel and discretization the resulting surface on a fixed grid. We follow the construction of \cite{images} and refer to that work for full methodological details. As noted in \cite{images}, the method is robust to the choice of resolution and Gaussian variance parameters. Hence, we've chosen suitable parameters and fixed the Gaussian variance at $0.1$ and use a $50 \times 50$ grid resolution throughout. No weighting function is applied, as we aim to evaluate the estimator without assuming higher-persistence features are more important. All persistence image summaries are computed using \texttt{GUDHI} \cite{GUDHI}.

%For both persistence diagram-- and persistence image--based representations, we also consider three variants: using only \(H_0\) persistence points, using only \(H_1\) persistence points, and using the combined set of \(H_0\) and \(H_1\) points. In the latter case, persistence points from different homological dimensions are pooled into a single representation. These variants are denoted in figures by \(\mathrm{PD}_{H_0}\), \(\mathrm{PD}_{H_1}\), and \(\mathrm{PD}\) for persistence diagram representations, and analogously by \(\mathrm{PI}_{H_0}\), \(\mathrm{PI}_{H_1}\), and \(\mathrm{PI}\) for persistence images.

\subsubsection{Results}
Predictive performance is evaluated using the coefficient of determination $R^2$ between the predicted and true parameter values. Let $n_{\mathrm{train}}$ and $n_{\mathrm{test}}$ denote the number of training and testing samples respectively. We fix \(n_{\mathrm{test}} = 50\) and vary \(n_{\mathrm{train}}\) to study performance under differing levels of data availability. For each value of \(n_{\mathrm{train}}\), the experiment is repeated \(10\) times with independent training and test sets, and we report the mean and standard deviation of the resulting \(R^2\) scores.

\begin{figure}[H]
  \centering
  \includegraphics[width=1\textwidth]{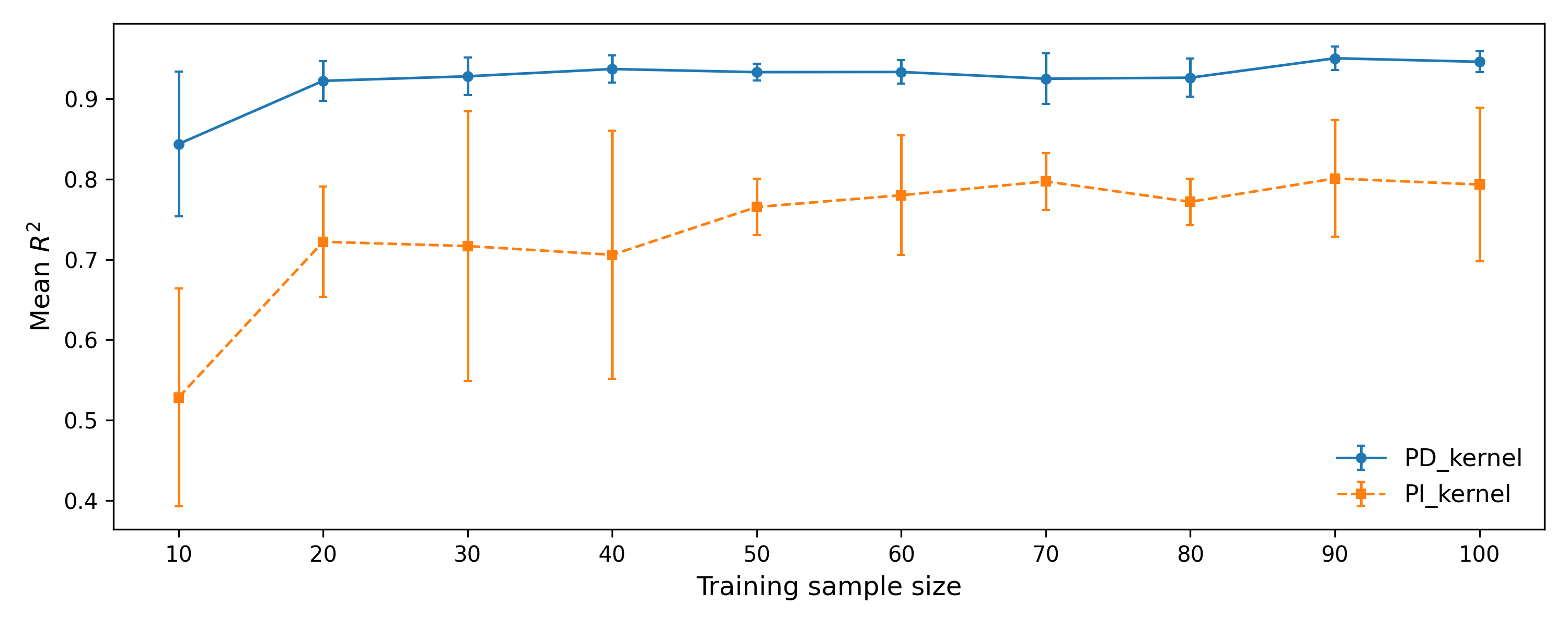}
  \caption{Mean predictive performance ($R^2$) as a function of training sample size $n_{\mathrm{train}}$ for persistence diagram (PD) and persistence image (PI) based methods. Points denote empirical means across repetitions. Error bars indicate ±1 standard deviation.}
  \label{fig:interpolation_results}
\end{figure}

Figure ~\ref{fig:interpolation_results} summarizes the results by plotting the mean \(R^2\) as a function of \(n_{\mathrm{train}}\) for each method, with error bars indicating one standard deviation across the \(10\) repetitions. Across all training sample sizes, the proposed kernel regression method on persistence diagrams attains higher mean \(R^2\) values than the persistence image based approach, with the performance gap most pronounced in low-sample regimes. As \(n_{\mathrm{train}}\) increases, the performance of both methods improves, though the persistence diagram based estimator remains consistently superior. The results indicate that preserving the geometry of persistence diagram space can yield improved predictive performance, particularly when training data are limited.

Exact numerical values corresponding to the curves in Figure \ref{fig:interpolation_results} are reported in the Table \ref{table_interpolation}.

\begin{comment}

\begin{table}[H]
\centering
\caption{Interpolation parameter experiment: mean $R^2$ with standard deviation (in parentheses) over 10 repetitions.}
\label{table_interpolation}
\begin{tabular}{ccc}
\toprule
\textbf{$n_{\mathrm{train}}$} 
& \textbf{PD\_kernel} 
& \textbf{PI\_kernel} \\
\midrule
10  & $0.8435\;(0.0901)$ & $0.5283\;(0.1354)$ \\
20  & $0.9221\;(0.0246)$ & $0.7219\;(0.0684)$ \\
30  & $0.9277\;(0.0233)$ & $0.7166\;(0.1675)$ \\
40  & $0.9368\;(0.0167)$ & $0.7058\;(0.1545)$ \\
50  & $0.9329\;(0.0101)$ & $0.7652\;(0.0351)$ \\
60  & $0.9331\;(0.0145)$ & $0.7798\;(0.0745)$ \\
70  & $0.9247\;(0.0317)$ & $0.7970\;(0.0353)$ \\
80  & $0.9260\;(0.0238)$ & $0.7716\;(0.0290)$ \\
90  & $0.9501\;(0.0146)$ & $0.8007\;(0.0724)$ \\
100 & $0.9458\;(0.0129)$ & $0.7931\;(0.0956)$ \\
\bottomrule
\end{tabular}
\end{table}

\begin{table}[H]
\centering
\caption{Curvature Experiment Results}
\label{table_curvature}
\begin{tabular}{lc}
\toprule
\textbf{Method} & \textbf{MSE} \\
\midrule
\texttt{ml\_on\_ph0}         & $0.3498$ \\
\texttt{ml\_on\_ph1}         & $0.6315$ \\
\texttt{ml\_on\_ph0\_longest}& $0.2912$ \\
\texttt{ml\_on\_ph1\_longest}& $0.6189$ \\
\texttt{ph0\_ml}             & $0.2669$ \\
\texttt{ph1\_ml}             & $0.5247$ \\
\texttt{ml\_on\_dis\_mat}    & $0.3916$ \\
\texttt{nn\_shallow}         & $1.1812$ \\
\texttt{nn\_deep}            & $2.7978$ \\
\texttt{point\_net}          & $37.1275$ \\
\bf{\texttt{PD\_kernel}}          & $\mathbf{0.2258}$ \\
\bottomrule
\end{tabular}
\end{table}

\end{comment}

\subsection{Curvature Parameter Experiment}
Next, we evaluate our method on a regression task introduced in \cite{PH_detects_Curvature}, which aims to predict the curvature $\kappa$ of an underlying surface with constant curvature from a point cloud sampled from a unit disk embedded in that surface. This task is motivated by applications in geometry processing and shape analysis, where curvature encodes fundamental local geometric information. An extensive computational study of this task is conducted in \cite{onPHeffectiveness}. We replicate their experimental setup, but reduce both the number of training samples and the number of points sampled per point cloud. This modification allows us to assess performance in a low-data regime, where differences between methods are most apparent. 

\subsubsection{Dataset}
Following the procedure of \cite{onPHeffectiveness}, each point cloud is generated by sampling from a unit disk $D_\kappa$ embedded in a surface with constant curvature $\kappa \in [-2, 2]$. For $50$ evenly spaced curvature values, we generate point clouds consisting of $50$ points, where the sampling density is proportional to the surface area measure. This yields $50$ training samples. Test data are constructed analogously using $50$ curvature values sampled uniformly at random from $[-2,2]$. We refer the reader to \cite{onPHeffectiveness} for full details of the data generation procedure and the associated code.

\subsubsection{Methods}
We evaluate the same kernel regression framework described in Section~\ref{implementation_steps}, using persistence diagram representations of the point clouds (denoted as \texttt{PD\_kernel}). For comparison, all baseline methods from \cite{onPHeffectiveness} are included. These are comprised of TDA-based machine learning approaches after vectorization (\texttt{ml\_on\_ph0}, \texttt{ml\_on\_ph1}, \texttt{ml\_on\_ph0\_longest}, \texttt{ml\_on\_ph1\_longest}, \texttt{ph0\_ml},
\texttt{ph1\_ml}, \texttt{ml\_on\_dis\_mat}) and neural network approaches (\texttt{nn\_shallow}, \texttt{nn\_deep}, \texttt{point\_net}). All methods are evaluated using the original implementations provided by the authors. We refer the reader to \cite{onPHeffectiveness} for full methodological details.

\subsubsection{Results}
Predictive performance is evaluated using the mean squared error (MSE), following the same evaluation protocol as \cite{onPHeffectiveness}. Table \ref{table_curvature} summarizes the results. The proposed kernel regression method on persistence diagrams attains the lowest MSE among the methods considered.

\begin{table}[H]
\centering

\begin{minipage}[t]{0.48\textwidth}
\centering
\captionof{table}{Interpolation parameter experiment: mean $R^2$ with standard deviation (in parentheses) over 10 repetitions.}
\label{table_interpolation}

\small
\begin{tabular}{ccc}
\toprule
\textbf{$n_{\mathrm{train}}$}
& \textbf{PD\_kernel}
& \textbf{PI\_kernel} \\
\midrule
10  & $0.8435\;(0.0901)$ & $0.5283\;(0.1354)$ \\
20  & $0.9221\;(0.0246)$ & $0.7219\;(0.0684)$ \\
30  & $0.9277\;(0.0233)$ & $0.7166\;(0.1675)$ \\
40  & $0.9368\;(0.0167)$ & $0.7058\;(0.1545)$ \\
50  & $0.9329\;(0.0101)$ & $0.7652\;(0.0351)$ \\
60  & $0.9331\;(0.0145)$ & $0.7798\;(0.0745)$ \\
70  & $0.9247\;(0.0317)$ & $0.7970\;(0.0353)$ \\
80  & $0.9260\;(0.0238)$ & $0.7716\;(0.0290)$ \\
90  & $0.9501\;(0.0146)$ & $0.8007\;(0.0724)$ \\
100 & $0.9458\;(0.0129)$ & $0.7931\;(0.0956)$ \\
\bottomrule
\end{tabular}
\end{minipage}
\hfill
\begin{minipage}[t]{0.48\textwidth}
\centering
\captionof{table}{Curvature experiment results: Mean Square Error (MSE) comparison between all methods.}
\label{table_curvature}

\small
\begin{tabular}{lc}
\toprule
\textbf{Method} & \textbf{MSE} \\
\midrule
\texttt{ml\_on\_ph0}          & $0.3498$ \\
\texttt{ml\_on\_ph1}          & $0.6315$ \\
\texttt{ml\_on\_ph0\_longest} & $0.2912$ \\
\texttt{ml\_on\_ph1\_longest} & $0.6189$ \\
\texttt{ph0\_ml}              & $0.2669$ \\
\texttt{ph1\_ml}              & $0.5247$ \\
\texttt{ml\_on\_dis\_mat}     & $0.3916$ \\
\texttt{nn\_shallow}          & $1.1812$ \\
\texttt{nn\_deep}             & $2.7978$ \\
\texttt{point\_net}           & $37.1275$ \\
\bf{\texttt{PD\_kernel}}      & $\mathbf{0.2258}$ \\
\bottomrule
\end{tabular}
\end{minipage}

\end{table}

\begin{comment}
    \subsection{Tecator Experiment} In this experiment we evaluate our kernel regression estimator on the \emph{Tecator} dataset which is a classical benchmark dataset used in functional data analysis.  

\subsubsection{Dataset} The Tecator dataset consists of $215$ meat samples, each accompanied by a near-infrared (NIR) absorbance spectrum. The objective is to predict the fat content of each sample, determined through laboratory measurements, from its corresponding spectral measurements.

To construct persistence diagrams, we treat each absorbance spectrum as a one-dimensional signal and compute its second derivative. We then apply a sublevel-set filtration to the resulting signal and compute its zero-dimensional persistent homology. Each observation is consequently represented by a persistence diagram encoding the multiscale structure of its spectral profile.

The resulting regression problem consists of predicting the fat content from the persistence diagram associated with each sample.

\subsubsection{Methods}

\subsubsection{Results}

\end{comment}

\subsection{Reynolds Number Experiment} 

Finally, we introduce a new experiment to investigate the applicability of the proposed method to regression problems involving complex spatial data. Specifically, we consider the prediction of the Reynolds number of a simulated fluid flow from its velocity field. The Reynolds number characterizes the relative importance of inertial and viscous forces in a fluid, with low Reynolds numbers generally associated with smooth, orderly flow patterns, while high Reynolds numbers give rise to more complex, irregular spatial structures and turbulent behaviour. This provides an interesting setting for persistence diagram regression, as the topological features of the velocity field may encode information about these underlying differences in flow structure.

\subsubsection{Dataset} 

We consider the 2D Unsteady Fluid Simulation Ensemble for Machine Learning introduced in \cite{Jakob2021FluidDataset}. The dataset contains $8000$ simulations of incompressible flow on a periodic square domain (a torus), each stored as a sequence of velocity field snapshots $(u,v)$ on a $512\times512$ grid at times $t\in[0,10]$. The simulations differ in their Reynolds number $\mathrm{Re} \in [1, 4096]$, which controls the turbulence and spatial complexity of the flow. Representative flows for different Reynolds numbers at different snapshot times are shown in Figure~\ref{fig:fluid_sample_data}. See \cite{Popinet2004Gerris} for details on how the fluid flow simulations were created. The regression task is to predict $\ln \mathrm{Re}$ of a simulation from a snapshot of the fluid velocity field at a fixed point in time for $t\in\{1,\,2.5,\,5,\,7.5,\,9\}$.

\subsubsection{Methods} 
\begin{figure}[H]
    \centering
    \includegraphics[width=1\textwidth]{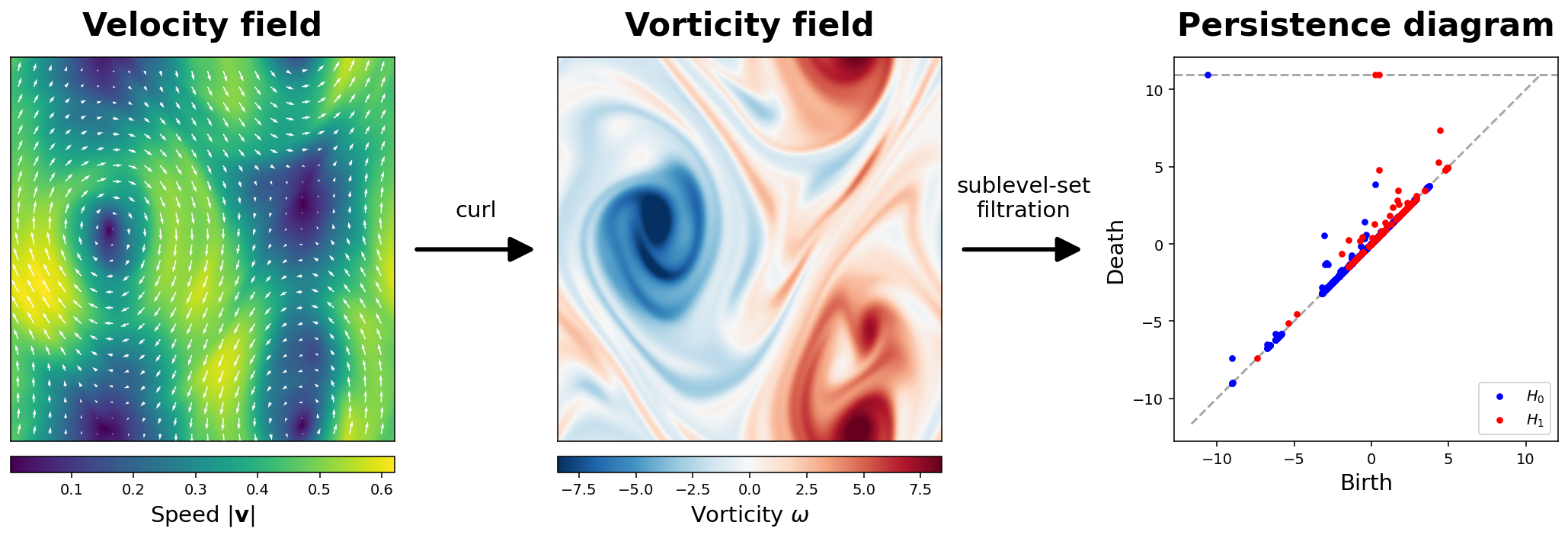}
    \caption{Data pipeline for converting a fluid simulation snapshot at a fixed time point to a persistence diagram. The velocity field (left) gives the vorticity field (middle) via its curl; and the sublevel-set filtration of the vorticity yields the persistence diagram (right).}
    \label{fig:fluid_data_pipeline}
\end{figure}

We use the same kernel regression procedure of Section~\ref{implementation_steps} with a few changes to adapt to this setting. Since the snapshots are vector fields on a grid, rather than point clouds, Step~1 differs from the other experiments. For each snapshot of the velocity field we compute the vorticity scalar field $\omega = \partial_x v - \partial_y u$ and the persistence diagram of its sublevel-set filtration (Figure~\ref{fig:fluid_data_pipeline}). Steps~2--4 in the procedure are unchanged.

The diagrams are computed with the cubical complex module of \texttt{GUDHI} \cite{GUDHI}. The $512\times512$ grid of vorticity values is passed as the filtration values on the vertices of a cubical complex, in which each higher-dimensional cell takes the maximum value of its vertices, so the sublevel set $\{\omega \le s\}$ is the union of the cells with filtration value at most $s$. Because the simulation domain is periodic, we use the periodic cubical complex with periodic boundary conditions in both directions, so the underlying space is a torus rather than a square.

\begin{figure}
    \centering
    \includegraphics[width=1\textwidth]{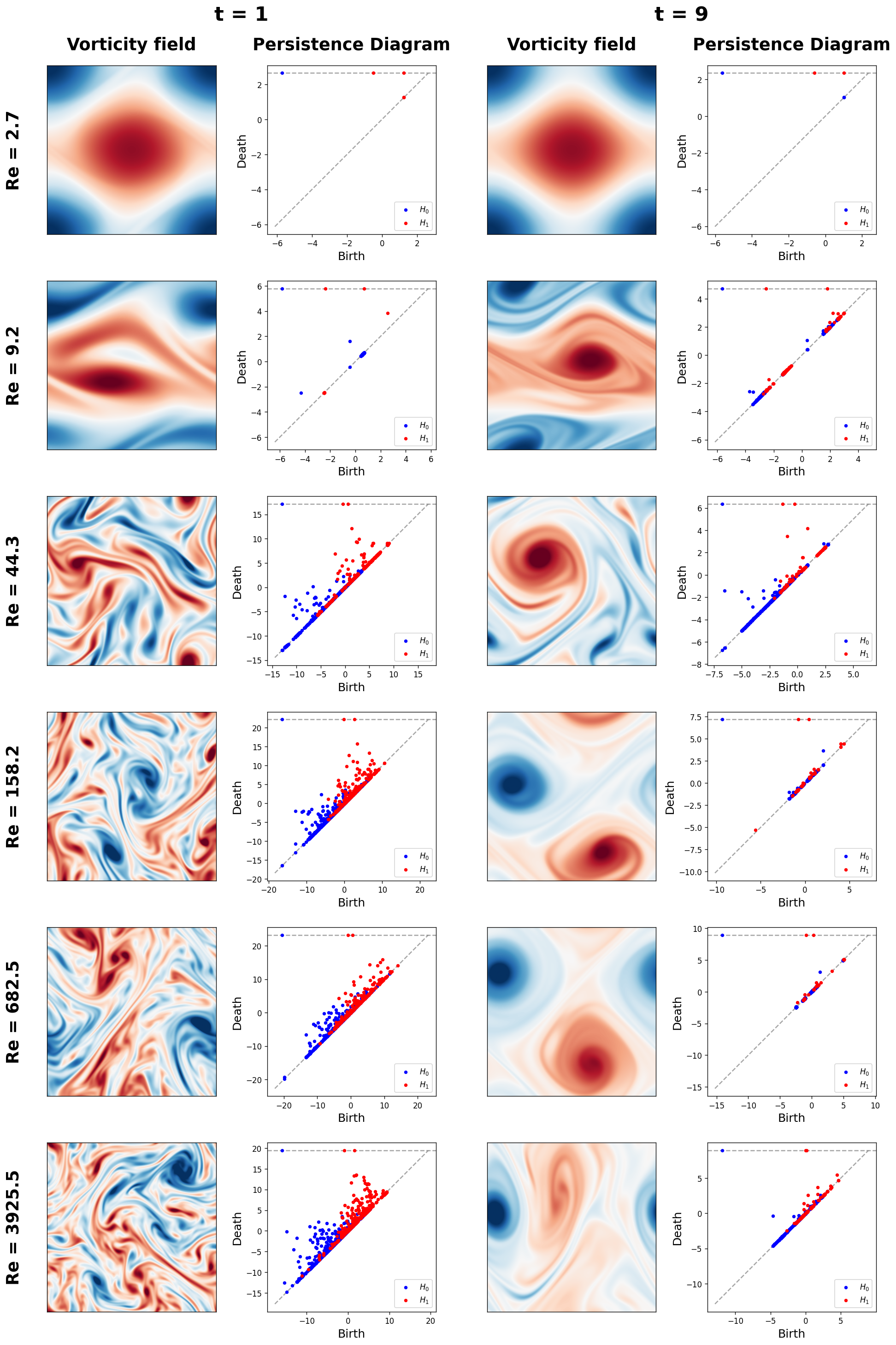}
    \caption{Sample vorticity fields from fluid dynamics dataset \cite{Jakob2021FluidDataset} and their associated persistence diagrams, at snapshot times $t=1$ and $t=9$. Different simulations are shown (with each row) with varying Reynolds numbers.}
    \label{fig:fluid_sample_data}
\end{figure}

Representative vorticity fields and persistence diagrams at different Reynolds numbers and snapshot times are shown in Figure~\ref{fig:fluid_sample_data}. 

\subsubsection{Results} We use $400$ simulations spread evenly over the whole range of Reynolds numbers, and the same simulations at every time point $t\in\{1,\,2.5,\,5,\,7.5,\,9\}$.  Performance is measured by the out-of-sample $R^2$ on $\ln \mathrm{Re}$. For each of $5$ random train/test splits, we hold out $100$ simulations for testing and train on nested subsets of the remaining simulations of size $n_{\mathrm{train}} \in \{25, 50, 100, 200, 300\}$.

\begin{figure}[H]
    \centering
    \includegraphics[width=1\textwidth]{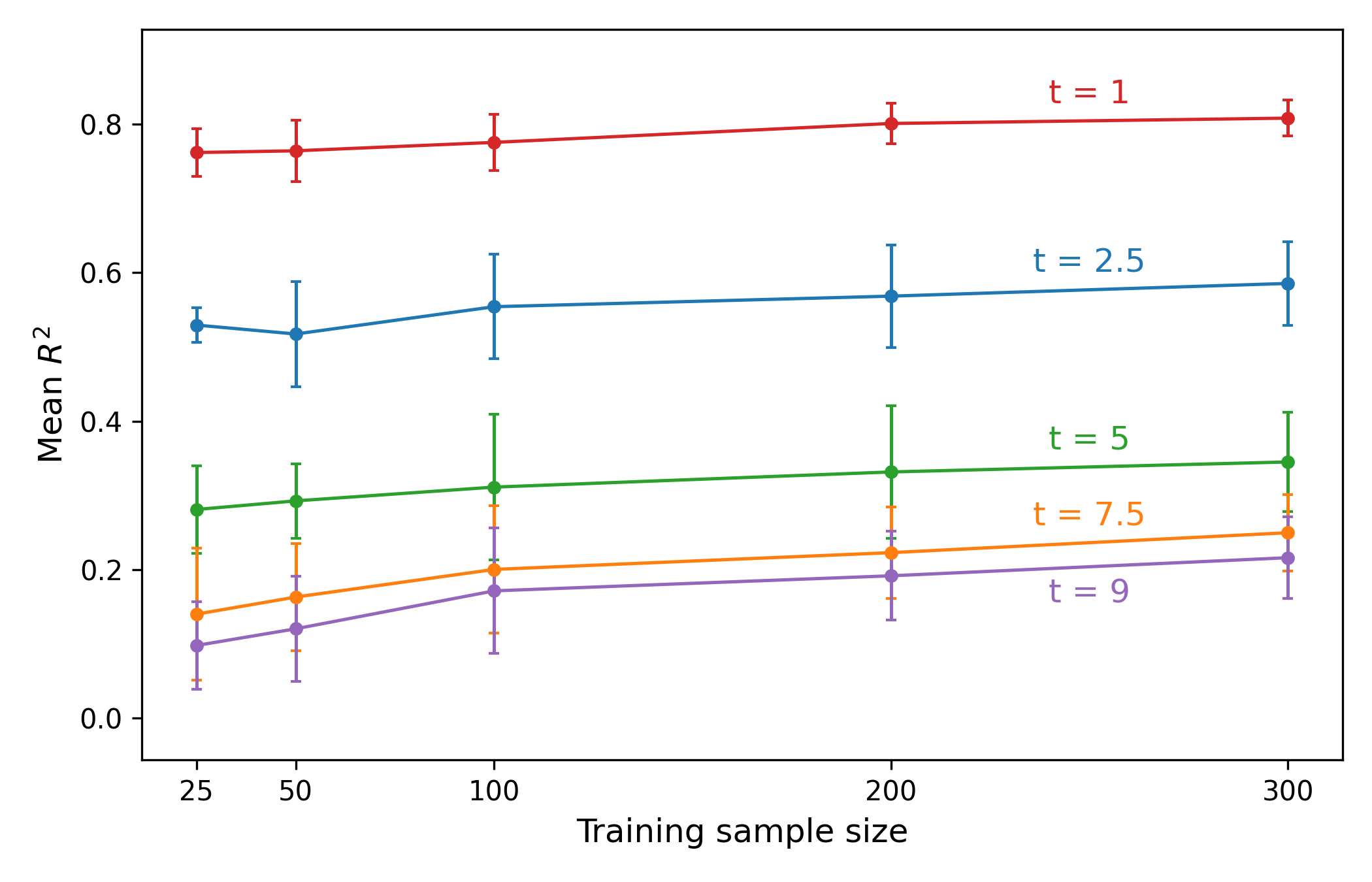}
    \caption{Mean predictive performance ($R^2$) as a function of training sample size $n_{\mathrm{train}}$ for persistence diagram kernel regression of $\ln \mathrm{Re}$ from a single vorticity snapshot, shown for snapshot times $t \in \{1, 2.5, 5, 7.5, 9\}$. Points denote empirical means across 5 random train/test splits with 100 test simulations each. Error bars indicate $\pm 1$ standard deviation.}
    \label{fig:results_fluid}
\end{figure}

Figure~\ref{fig:results_fluid} plots the mean $R^2$ as a function of $n_{\mathrm{train}}$ for each snapshot time, with error bars indicating one standard deviation across the $5$ splits, and the exact values are reported in Table~\ref{table_fluid}. The mean $R^2$ is higher at earlier snapshot times for every training size, demonstrating that our kernel method for persistence diagrams has strong predictive power for inferring the Reynold's number at time $t=1$.

For snapshots of the fluid simulations at later times, the method does not perform as well, and the $R^2$ value drops significantly. On further examination, we see that this is due to the fluid settling in to an almost-steady state at later times with similar vortex structures in each simulation at different Reynold's numbers. See Figure~\ref{fig:fluid_sample_data} to compare the vorticity fields at $t=1$ and $t=9$. 

\begin{figure}[H]
    \centering
    \includegraphics[width=1\textwidth]{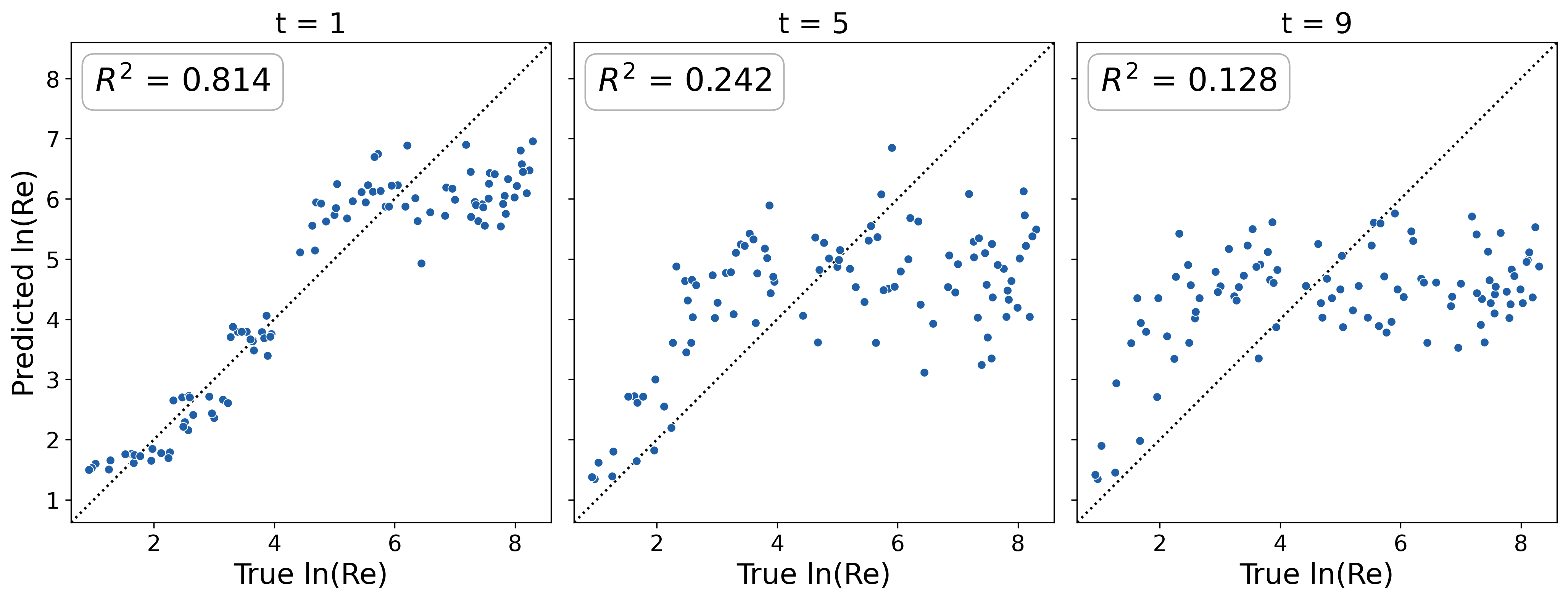}
    \caption{Predicted versus true $\ln \mathrm{Re}$ on a single train/test split ($n_{\mathrm{train}}=300$, $n_{\mathrm{test}}=100$) for persistence diagram kernel regression at snapshot times $t=1$, $5$ and $9$. Each point is one test simulation, and the dotted line $y=x$ marks perfect prediction. The $R^2$ of the split is shown in each panel.}
    \label{fig:results_fluid_scatter}
\end{figure}

The similarity at say $t=9$ makes the corresponding persistence diagrams harder to distinguish, but that difficulty also holds to some extent for visual inspection of the vortex field. Looking at the vortex fields in Figure~\ref{fig:fluid_sample_data}, it is hard to see much difference for example between the $\mathrm{Re}=158.2$ and $\mathrm{Re}=682.5$ simulations at $t=9$. It is no surprise then that our kernel estimation of the Reynold's number loses accuracy at later times.
%A possible explanation is that the flow is more chaotic at early times, so that the persistence diagram of the vorticity field captures rich structure that depends on the Reynolds number. At later times the flow settles into similar vortex structures which makes the corresponding persistence diagrams harder to distinguish. 
%See Figure~\ref{fig:fluid_sample_data} to compare the vorticity fields at $t=1$ and $t=9$, which also shows that the vortex structures look more similar to each other at higher Reynolds numbers. 
This is reflected in Figure~\ref{fig:results_fluid_scatter}, which plots the predicted against the true $\ln \mathrm{Re}$ for a single split with $n_{\mathrm{train}}=300$ at $t=1$, $5$ and $9$. The predicted values deviate much more from the true values at higher Reynolds numbers, and this deviation is larger at the later time points ($R^2=0.814$, $0.242$ and $0.128$ respectively).

\begin{table}[H]
\centering
\caption{Reynolds Number Experiment: mean $R^2$ with standard deviation (in parentheses) over 5 random train/test splits, for each snapshot time $t$.}
\label{table_fluid}
\begin{tabular}{cccccc}
\toprule
\textbf{$n_{\mathrm{train}}$}
& \textbf{$t=1$}
& \textbf{$t=2.5$}
& \textbf{$t=5$}
& \textbf{$t=7.5$}
& \textbf{$t=9$} \\
\midrule
25  & $0.7620\;(0.0319)$ & $0.5295\;(0.0232)$ & $0.2811\;(0.0590)$ & $0.1401\;(0.0886)$ & $0.0979\;(0.0589)$ \\
50  & $0.7643\;(0.0411)$ & $0.5175\;(0.0710)$ & $0.2927\;(0.0502)$ & $0.1632\;(0.0720)$ & $0.1204\;(0.0708)$ \\
100 & $0.7756\;(0.0382)$ & $0.5543\;(0.0704)$ & $0.3113\;(0.0980)$ & $0.2004\;(0.0856)$ & $0.1715\;(0.0844)$ \\
200 & $0.8011\;(0.0275)$ & $0.5685\;(0.0689)$ & $0.3317\;(0.0893)$ & $0.2230\;(0.0619)$ & $0.1918\;(0.0597)$ \\
300 & $0.8083\;(0.0243)$ & $0.5856\;(0.0563)$ & $0.3451\;(0.0669)$ & $0.2498\;(0.0515)$ & $0.2162\;(0.0552)$ \\
\bottomrule
\end{tabular}
\end{table}

\clearpage
\appendix

\section{Technical Lemmas}
The lemma below allows us to bound probabilities with covering numbers of our space. The original proof can be found in \cite{Abraham}.
\begin{lemma}[Small Ball Probability Control]
    Let $(M,\mathcal{B}(M), \mathbb{P})$ be a probability space and let $S^* \subset M$ be compact. Consider the ball $B(x,r) \subset M$ for some fixed $x \in M$ and $r > 0$. Then:
    \begin{align*}
        \int_{S^*} \frac{1}{\mathbb{P}(B(x, r))} d\mathbb{P} \leq \mathcal{C}\left(S^*, \frac{r}{2}\right) 
    \end{align*}
    \label{SBP_control_compact}
\end{lemma}

The proof for the lemma below can be found in \cite{Gyofri_book}.
\begin{lemma}[Binomial inverse expectation bound]\label{lemma:binomial_inverse}
Let $B(n,p)$ be a binomial random variable with parameters $n$ and $p$. Then:
\begin{enumerate}[(i)]
    \item 
    \[
        \mathbb{E}\!\left[\frac{1}{1 + B(n,p)}\right] \le \frac{1}{(n+1)p}
    \]
    \item 
    \[
        \mathbb{E}\!\left[\frac{1}{B(n,p)} \mathds{1}_{\{B(n,p) > 0\}}\right] \le \frac{2}{(n+1)p}
    \]
\end{enumerate}
\end{lemma}

\section{Additional Proofs}

\emph{Proof of Theorem \ref{thm:Stones}}
\begin{proof}
\label{proof_for_stones}
\begin{align*}
    \lVert \hat{m}_n(X) - m(X) \rVert_{L^2}^2 &= \mathbb{E}\left[|\hat{m}_n(X) - m(X)|^2\right] \\
    &= \mathbb{E}\left[\left(\sum_{i=1}^n W_i(X)Y_i - m(X)\right)^2\right]
\end{align*}
We can write the above expression as:
\begin{align*}
&\mathbb{E}\left[\left(\sum_{i=1}^n W_i(X)Y_i - \sum_{i=1}^n W_i(X)m(X_i) + \sum_{i=1}^n W_i(X)m(X_i) \right. \right. \\
& \left. \left. \quad - \sum_{i=1}^n W_i(X)m(X) + \sum_{i=1}^n W_i(X)m(X) - m(X)\right)^2\right]
\end{align*}
Using the inequality \((a + b + c)^2 \leq 3a^2 + 3b^2 + 3c^2\), we split the expectation as follows:
\begin{align*}
    \mathbb{E}\left[(\hat{m}_n(X) - m(X))^2\right] \leq &\, \, \, 3\mathbb{E}\left[\left(\sum_{i=1}^n W_i(X)(Y_i - m(X_i))\right)^2\right] \\
    &+ 3\mathbb{E}\left[\left(\sum_{i=1}^n W_i(X)(m(X_i) - m(X))\right)^2\right] \\
    &+ 3\mathbb{E}\left[\left(\sum_{i=1}^n W_i(X) - 1\right)^2 m(X)^2\right] \quad \\
    &=: 3I_n + 3J_n + 3L_n
\end{align*}
We will show that each of \(I_n, J_n, L_n \to 0\) as \(n \to \infty\), which is the desired result.
\begin{enumerate}
    \item First we show \(L_n \to 0\) as \(n \to \infty\). We know by condition 1, that almost surely for all \(n\):
    \begin{align*}
         \left( \sum_{i=1}^n W_i(X) - 1 \right)^2 m(X)^2  &\leq (D+1)^2 |m(X)|^2
    \end{align*}
    We only need to consider $m$ on the compact support $S^*$. Since $m$ is continuous on $S^*$, we know $|m(X)| \leq B_m$ for some $B_m \geq 0$. Hence, we attain a uniform bound by an integrable function:
    \begin{align*}
        \left( \sum_{i=1}^n W_i(X) - 1 \right)^2 m(X)^2 \leq (D+1)^2B_m^2
    \end{align*}
    Finally, since condition 3 gives us the following convergence in probability:
    \begin{align*}
        \left( \sum_{i=1}^n W_i(X) - 1 \right)^2 m(X)^2 \xrightarrow{\mathbb{P}} 0 \quad \text{as} \quad n \to \infty
    \end{align*}
    we can apply the Vitali convergence theorem to obtain the desired result:
    \begin{align*}
        L_n &= \mathbb{E}\left[\left(\sum_{i=1}^n W_i(X) - 1\right)^2 m(X)^2\right] \to 0 \\
    \end{align*}

\item  Next we work on \(I_n \to 0\) as \(n \to \infty\)
    \begin{align*}
        I_n &= \mathbb{E}\left[\left(\sum_{i=1}^n W_i(X)(Y_i - m(X_i))\right)^2\right] \\
        &= \sum_{i=1}^n \sum_{j=1}^n \mathbb{E}\left[\, W_i(X)W_j(X)(Y_i - m(X_i))(Y_j - m(X_j)) \, \right]
    \end{align*}
But the cross-terms disappear since for \(i \neq j\):
\begin{align*}
    &\mathbb{E} \left[ W_i(X)\cdot W_j(X)\cdot (Y_i - m(X_i))\cdot (Y_j - m(X_j)) \right] \\
    &\hspace{3em}= \mathbb{E} \left[ \, \mathbb{E}\left[ (Y_j - m(X_j)) \mid X, X_1, \dots, X_n, Y_i \right] \, W_i(X) W_j(X) (Y_i - m(X_i)) \, \right] \\
    & \hspace{12em} \text{by law of total expectation and independence} \\
    &\hspace{3em}= 0 \quad \text{since} \quad \mathbb{E} \left[ (Y_j - m(X_j)) \mid X, X_1, \dots, X_n, Y_i \right] = \mathbb{E}(Y_j \mid X_j) - m(X_j) = 0.
\end{align*}
So it remains to control:
\begin{align*}
    I_n &= \mathbb{E}\left[\sum_{i=1}^n W_i(X)^2(Y_i - m(X_i))^2\right] = \mathbb{E}\left[\sum_{i=1}^n W_i(X)^2 \epsilon_i^2 \right]
\end{align*}
But conditioning on $X, X_1, ..., X_n$ gives:
\begin{align*}
    \mathbb{E}\left[\sum_{i=1}^n W_i(X)^2\epsilon_i^2\right] &= \mathbb{E}\left[\, \sum_{i=1}^n \mathbb{E}\left( W_i(X)^2\epsilon_i^2 \mid X, X_1, ..., X_n\right)\right] \\
    &= \mathbb{E} \left[ \, \sum_{i=1}^n W_i(X)^2 \, \mathbb{E}\left( \epsilon_i^2 \mid X_i\right)  \, \right] \\
    &\leq \sigma^2 \, \mathbb{E} \left[ \, \sum_{i=1}^n W_i(X)^2 \,  \right] \xrightarrow{n \to \infty} 0
\end{align*}
where the inequality follows from the conditional
second-moment assumption in the regression set up, and the convergence follows from condition 4. Hence $I_n \to 0$ as required. \\

\item Now it remains to show $J_n \to 0$ 
\begin{align*}
    J_n &\leq \mathbb{E}\left[\left(\sum_{i=1}^n |W_i(X)| \cdot |m(X_i) - m(X)|\right)^2\right] \\
    &= \mathbb{E}\left[\left(\sum_{i=1}^n |W_i(X)|^{\frac{1}{2}} \cdot |W_i(X)|^{\frac{1}{2}} \, |m(X_i) - m(X)|\right)^2\right]
\end{align*}
By the Cauchy-Schwarz inequality and applying condition 1:
\begin{align*}
    &\leq \mathbb{E}\left[\left( \sum_{i=1}^n |W_i(X)| \right) \left( \sum_{i=1}^n |W_i(X)| \cdot |m(X_i) - m(X)|^2 \right)\right] \\
    &\leq D \cdot \mathbb{E}\left[\sum_{i=1}^n |W_i(X)| \cdot |m(X_i) - m(X)|^2\right]
\end{align*}
We now split the expectation depending on whether \(X_i\) lies inside or outside the ball \(B(X, \delta)\):
\begin{align*}
    &D \cdot \mathbb{E}\left[\sum_{i=1}^n |W_i(X)| \cdot |m(X_i) - m(X)|^2 \cdot \mathds{1}_{\{d(X_i, X) > \delta\}}\right] \\
    &\hspace{5em} + D \cdot \mathbb{E}\left[\sum_{i=1}^n |W_i(X)| \cdot |m(X_i) - m(X)|^2 \cdot \mathds{1}_{\{d(X_i, X) \leq \delta\}}\right]
\end{align*}
The first expectation goes to zero by application of condition 2 and since $|\, m(X) \,| \leq b$:
\begin{align*}
    &D\cdot \mathbb{E}\left[\sum_{i=1}^n |W_i(X)| \cdot |m(X_i) - m(X)|^2 \cdot \mathds{1}_{\{d(X_i, X) > \delta\}}\right] \\
    &\hspace{5em} \leq 4b^2 D \cdot \mathbb{E}\left[\sum_{i=1}^n |W_i(X)| \cdot \mathds{1}_{\{d(X_i, X) > \delta\}}\right] \xrightarrow{n \to \infty} 0
\end{align*}
The second expression (for values inside the ball $B(X,\delta)$) can be controlled because:
\begin{align*}
    &D \cdot \mathbb{E}\left[\sum_{i=1}^n |W_i(X)| \cdot |m(X_i) - m(X)|^2 \cdot \mathds{1}_{\{d(X_i, X) \leq \delta\}}\right] \\
    &\hspace{5em} \leq D^2 \sup_{\{p_1, p_2 \in \, S^* \, : \, d(p_1, p_2) \leq \delta\}} \Bigl\{ |m(p_1) - m(p_2)|^2 \Bigr\} 
\end{align*}
And since $m$ is continuous on a compact set, we know it is uniformly continuous. Thus:
\begin{align*}
    \sup_{\{p_1, p_2 \in \, S^* \, : \, d(p_1, p_2) \leq \delta\}} \Bigl\{ |m(p_1) - m(p_2)|^2 \Bigr\} \to 0 \quad \text{as } \delta \to 0
\end{align*}
And so $J_n \xrightarrow{n \to \infty} 0$ as required.
\end{enumerate}

\noindent Thus, we have shown that \(I_n, J_n, L_n \to 0\), completing the proof.
\end{proof}

\emph{Proof for Theorem \ref{consistency_theorem}}
\begin{proof}
\label{proof_for_consistency}
We just need to verify that the weights for the kernel estimator satisfy Stone's conditions. Note that we use the following notation for scaling the kernel function:
\begin{align*}
    K_{h}(x) = K\left(\frac{x}{h}\right)
\end{align*}
\begin{enumerate}

\item Since our weights are normalised in the definition this condition becomes trivial since for any $n$ and $h_n$:
    \begin{align*}
    \sum_{i = 1}^{n}|W_i(X)| = \frac{ \sum_{i = 1}^{n} K_{h_n}(d(X_i, X))}{\sum_{j=1}^n K_{h_n}(d(X_j, X))}\mathds{1}_{A} = \mathds{1}_{A} \leq 1 
    \end{align*}

\item By our condition on the kernel function we know:
\begin{align*}
    0 \leq K\left(\frac{d(X_i, X)}{h_n}\right) &\leq \mathds{1}_{\left\{\frac{d(X_i, X)}{h_n} \in [0, R]\right\}} \\
    &= \mathds{1}_{\left\{d(X_i, X) \in [0, Rh_n]\right\}}
\end{align*}
But since $h_n \xrightarrow{n \to \infty} 0$, for any $\delta > 0$ we can choose $n$ large enough so $Rh_n < \delta$ so:
\begin{align*}
    \mathds{1}_{\left\{d(X_i, X) > \delta \right\}} =  1 &\implies \mathds{1}_{\left\{d(X_i, X) \in [0, Rh_n]\right\}} = 0 \\
    &\implies K_{h_n}\left(d(X_i, X)\right) = 0 \\
    &\implies |W_i(X)| \mathds{1}_{\left\{d(X_i, X) > \delta \right\}} \xrightarrow{n \to \infty} 0
\end{align*}
Finally by dominated convergence theorem we obtain the desired result:
\begin{align*}
    \lim_{n \to \infty} \mathbb{E}\left[\sum_{i = 1}^{n} \left|W_i(X)\right| \cdot \mathds{1}_{\left\{d(X_i, X) > \delta \right\}} \right] &= \mathbb{E}\left[\lim_{n \to \infty} \sum_{i = 1}^{n} \left|W_i(X)\right| \cdot \mathds{1}_{\left\{d(X_i, X) > \delta \right\}} \right] \\
    &= 0
\end{align*}

\item Since our weights are normalised, they always sum to $1$ as long as the sum of kernel values is non-zero. We will show that the probability of the kernel sum being $0$ is also $0$ in the limit:
\begin{align*}
    \mathbb{P}\left(\sum_{i = 1}^{n} \left|W_i(X)\right| \neq 1\right) 
    &= \mathbb{P}\left(\sum_{i = 1}^{n} K_{h_n}(d(X_i, X)) = 0\right) \\
    &\leq \mathbb{P}\left(\sum_{i = 1}^{n} \mathds{1}_{\{X_i \in B(X, rh_n)\}} = 0\right)
\end{align*}
For any fixed $X = x$, we are demanding that there are no samples $X_i$ in $B(x, rh_n)$. Since the sample observations are i.i.d:
\begin{align*}
    \mathbb{P}\left(\sum_{i = 1}^{n} \mathds{1}_{\{X_i \in B(X, rh_n)\}} = 0\right) &= \int_{S^*} \Bigl( 1 - \mathbb{P}\bigl(B(x,rh_n)\bigr) \Bigr)^n \hspace{1em} d\mathbb{P} \\
    &\leq \int_{S^*} \text{exp}\Bigl(-n \cdot \mathbb{P}\bigl(B(x,rh_n)\bigr) \Bigr) \, d\mathbb{P} \\
    & \hspace{6em} \text{since} \, \, \forall x \in [0,1]: (1-x)^n \leq e^{-nx} \\
    & = \int_{S^*} \frac{n \cdot \mathbb{P}\bigl(B(x,rh_n)\bigr)}{n \cdot \mathbb{P}\bigl(B(x,rh_n)\bigr)} \text{exp}\Bigl(-n \cdot \mathbb{P}\bigl(B(x,rh_n)\bigr) \Bigr) \, d\mathbb{P} \\
    &\leq \int_{S^*} \frac{1}{n \cdot \mathbb{P}\bigl(B(x,rh_n)\bigr)} \, d\mathbb{P} \\
    &\hspace{6em} \text{since} \, \, \forall x \in [0,1]: nxe^{-nx} \leq 1
\end{align*}

Finally, we can apply the small ball probability control Lemma \ref{SBP_control_compact} and obtain:
\begin{align*}
    \frac{1}{n} \int_{S^*} \frac{1}{\mathbb{P}\bigl(B(x,rh_n)\bigr)} \, d\mathbb{P} \leq \frac{\mathcal{C}\left(S^*, \frac{rh_n}{2}\right)}{n} \xrightarrow{n \to \infty} 0
\end{align*}

\item Since $K(x) \leq 1$ we have:
\begin{align*}
    \sum_{i=1}^n W_i(X)^2 = \frac{ \sum_{i=1}^n \left(K_h(d(X_i, X))^2 \cdot \mathds{1}^2_{A}\right)}{\left( \sum_{j=1}^n K_h(d(X_j, X)) \right)^2} \leq \frac{\mathds{1}_A}{\sum_{i=1}^n K_h(d(X_i, X))}
\end{align*}
Since the observations are i.i.d for a fixed $X = x$, let $Z_n(x) = \sum_{i=1}^n \mathds{1}_{X_i \in B(x, rh_n)}$ be a binomial random variable where $p_n(x) = \mathbb{P}(B(x, rh_n))$ is the probability of success. 

When $Z_n(X) > 0$ we have:
\begin{align*}
    \sum_{i=1}^n W_i(X)^2 \leq \frac{1}{bZ_n(X)}
\end{align*}
since $b \cdot \mathds{1}_{[-r, r]}(x) \leq K(x)$. When $Z_n(X) = 0$:
\begin{align*}
    \sum_{i=1}^n W_i(X)^2 \leq \left(\sum_{i=1}^n W_i(X)\right)^2 \leq 1
\end{align*}
since the weights are non-negative and sum to at most $1$. Combining these two cases gives:
\begin{align*}
    \sum_{i=1}^n W_i(X)^2
    \leq
    \mathds{1}_{\{Z_n(X)=0\}}
    +\frac{1}{b}
    \frac{\mathds{1}_{\{Z_n(X)>0\}}}{Z_n(X)}.
\end{align*}
Taking the conditional expectation given $X = x$ and integrating over the distribution of $X$ for the expectation, we obtain:
\begin{align*}
    \mathbb{E}\left[
    \sum_{i=1}^n W_i(X)^2
    \right]
    &\leq
    \int_{S^*}
    \left[
    \mathbb{P}(Z_n(x)=0)
    +\frac{1}{b}
    \mathbb{E}\left(
    \frac{\mathds{1}_{\{Z_n(x)>0\}}}{Z_n(x)}
    \right)
    \right]d\mathbb{P}\\
    &=
    \int_{S^*}
    \left[
    (1-p_n(x))^n
    +\frac{1}{b}
    \mathbb{E}\left(
    \frac{\mathds{1}_{\{Z_n(x)>0\}}}{Z_n(x)}
    \right)
    \right]d\mathbb{P}\\
    &\leq
    \int_{S^*}
    \left[
    \frac{1}{np_n(x)}
    +\frac{2}{b(n+1)p_n(x)}
    \right]d\mathbb{P} \quad \text{by Lemma \ref{lemma:binomial_inverse}}\\
    &=\left(
    \frac{1}{n}+\frac{2}{b(n+1)}
    \right)
    \int_{S^*}\frac{1}{p_n(x)}\,d\mathbb{P}
\end{align*}
Again applying the small ball probability control Lemma \ref{SBP_control_compact} we obtain the desired result:
\begin{align*}
    \left(\frac{1}{n}+\frac{2}{b(n+1)}\right) \int_{S^*} \frac{1}{p_n(x)} \, d\mathbb{P} \leq \left(1+\frac{2}{b}\right) \cdot \frac{\mathcal{C}\left(S^*, \frac{rh_n}{2}\right)}{n} \xrightarrow{n \to \infty} 0
\end{align*}

\begin{comment}
    \begin{align*}
    \mathbb{E}\left[ \sum_{i=1}^n W_i(X)^2 \right] &\leq \frac{1}{b} \mathbb{E} \left[ \int_{S^*}\frac{\mathds{1}_A}{Z(x)} \, d\mathbb{P} \, \right] \\
    &\leq \frac{1}{b}  \int_{S^*} \mathbb{E} \left[ \frac{\mathds{1}_A}{Z(x)} \right] d\mathbb{P} \quad \text{by Fubini-Tonelli} \\
    &\leq \frac{1}{b}  \int_{S^*} \frac{2}{(n+1)p(x)} \, d\mathbb{P} \quad \text{by Lemma \ref{lemma:binomial_inverse}} \\
    & = \frac{2}{(n+1)b} \int_{S^*} \frac{1}{\mathbb{P}(B(x,rh_n))} \, d\mathbb{P}
    \end{align*}
\end{comment}

\end{enumerate}
\end{proof}

\bibliographystyle{imsart-number}
\bibliography{refs}

\end{document}